\documentclass{svjour3}
\smartqed

\usepackage{psfig}
\usepackage{amssymb}
\usepackage[usenames]{color}
\usepackage{url} 
\usepackage{soul} 

\journalname{Probability Theory and Related Fields}

\begin{document}

\title{The Reverse of the Law of Large Numbers}

\author{Kieran Kelly \and Przemys{\l}aw Repetowicz \and Seosamh macR{\'e}amoinn}

\institute{K. Kelly: \at 
Probability Dynamics, IFSC House, Custom House Quay, Dublin 1, Ireland.\\ \email{kkelly@probabilitydynamics.com} 
\and
P. Repetowicz: \at \email{prepetowicz@probabilitydynamics.com}
\and
S. macR{\'e}amoinn: \at \email{smacreamoinn@probabilitydynamics.com}
}

\date{Received: \today / Accepted: date}

\maketitle

\begin{abstract}
The Law of Large Numbers tells us that as the sample size ($N$) is increased, the sample mean converges on the population mean, provided that the latter exists. In this paper, we investigate the opposite effect: keeping the sample size fixed while increasing the number of outcomes ($M$) available to a discrete random variable. We establish sufficient conditions for the variance of the sample mean to increase monotonically with the number of outcomes, such that the sample mean \textit{``diverges''} from the population mean, acting like an \textit{``reverse''} to the law of large numbers. These results, we believe, are relevant to many situations which require \textit{sampling} of statistics of certain finite discrete random variables.


\end{abstract}

\section{Introduction}

In probability theory it is customary to investigate two broad families of problems, (i) the convergence of sums of large numbers of random variables, and (ii) the estimation of likelihoods of deviations of those sums from their large number limits \cite{Feller}.
The probability spaces of the random variables in question are usually fixed
and one investigates how fast do properties of a sample drawn from the statistical population \textit{\textbf{converge}} towards properties of that statistical population (to be called population for brevity, later on).

A good example of this is one of the most fundamental theorems of probability, the \textit{\textbf{Law of Large Numbers}} (LLN). The weak version can be stated (see~\cite{Grinstead} for example) as :
\begin{quotation}
Let $X_1, X_2, \ldots, X_N$ be an independent trials process, with finite expected value $\mu = E[X_i]$ and variance $\sigma^2 = Var(X_i)$. Let $S_N = \sum_{i=1}^{N} X_i$. Then $\forall \mbox{ } \epsilon > 0$,
\begin{equation}
\lim_{N \to \infty} P\left[ \left| \frac{S_N}{N} - \mu \right| \geq \epsilon \right] = 0 \label{eq:LLN}
\end{equation}
\end{quotation}
Noting that the quantity $\frac{S_N}{N}$ is nothing more than the sample mean, $\overline{X}_N$, the Law of Large Numbers can be stated in words as ``for independent trials of a sample of a distribution with finite population mean
, the sample mean should approach the population mean, as the number of trials (i.e. size of the sample) gets very large''.

In this paper we pose questions in a different way. We consider a sequence of populations from discrete, finite probability spaces whose size (to be termed number of outcomes in what follows) increases. We draw a sample of a fixed size from each of these populations, and analyze under what conditions, and how quickly, do the properties of the samples \textit{\textbf{diverge}} from properties of the populations.
In effect we are increasing monotonically the number of outcomes available to each trial, while keeping the number of trials, of sample size ($N$) constant.
In other words we formulate the \textbf{\textit{\textbf{``Reverse'' of the Law of Large Numbers}}}.

We illustrate this using a thought experiment based on a real life situation - betting on horses. 
Each horse is given a unique label.
Here the random variable in question is the label of the winning horse, the number of outcomes $M$ is the number of horses in the race, and the sample size $N$ stands for the number of repetitions of the race.
If the number of horses in the race is fixed ($M=\mbox{const}$) then, if the number of race repetitions becomes large ($N\rightarrow \infty$), the average label of the winning horse converges towards that of the expected value of the label of the winning horse. We have made the idealized assumption that the probability distribution itself remains fixed thus ignoring aging effects, weather changes, different tracks etc.
Now imagine that, as time progresses, new horses are added to the race indefinitely ($M\rightarrow \infty$) and that each time a new horse has been added, the same number of races ($N=\mbox{const}$) is held.
In this situation we intuitively perceive that the average label of the winning horse becomes less likely to be linked to the expected value of the label of the winning horse. 
In other words, the variance of the average label of the winner will increase at a rate that depends on the number of horses that participate in the race and on the probability distribution that a given horse wins.

Another thought experiment to illustrate this idea is to think of an unbiased $M$-sided ``die''. Here we ask what value appears face up when we throw the die. If we think of the standard die, with $M=6$, and consider $N = 1000$ throws, we would expect the average throw to be pretty close to the expected value of, in this case, $3.5$. However, if we increase $M$ to larger and larger values, but keep $N$ fixed, we would no longer expect to always get sample means close to the population mean (expected value).

In the above examples, we looked at one particular statistic - the sample mean. The Law of Large Numbers, in equation (\ref{eq:LLN}), tells us that the sample mean approaches the population mean as the number of trials increases. This is normally proven using Chebyshev's inequality (assuming that the variance, $\sigma^2$, is finite), which for the sample mean, is given by:
\begin{equation}
P\left[ \left| \overline{X}_N - \mu \right| \geq \epsilon \right] \leq \frac{\sigma^2}{N \epsilon^2} \label{eq:Cheby}
\end{equation} 
Here we used the fact that the variance of the sample mean is given by: $\sigma^2_{\overline{X}_N} = \frac{\sigma^2}{N}$. 
Chebyshev's inequality offers a sense of ``how large'' $N$ must be to see the desired convergence described by $\epsilon$. For fixed $N$ and $\epsilon$, the variance of the random variable, $X$, controls the bound on the probability of the sample mean being close to the population mean\footnote{Of course, this probability can never be greater than $1$ so bounds exceeding $1$ are ``loose''.}.

In this paper, we look at probability distributions of finite discrete random variables. We consider what happens when additional outcomes are deemed possible, but the distribution still retains the same basic form. We then have, in essence, a sequence consisting of probability distributions that are all similar, but with higher terms in the sequence corresponding to those distributions with a higher number of outcomes ($M$). We also have corresponding sequences for the expected values, and variances of these distributions. We are interested in those such sequences of variances which diverge as the number of outcomes becomes large. Then, from equation (\ref{eq:Cheby}), and for fixed $N$ and $\epsilon$, it is clear that sample means for the corresponding random variables would not be likely to be close to the expected values of those random variables.

Therefore, we seek conditions for these sequences of variances to diverge as the number of outcomes becomes very large. We also attempt to determine classes of (sequences of) distributions that correspond to particular rates of divergence of the variance. In the course of this work we will use the terms probability measure and probability distribution interchangeably. We will call the sequences of distributions described above, as simply distributions (that are functions of $M$), and the sequences of variances as variances (that are functions of $M$). We will attempt to formulate our considerations and results in the axiomatic language of Kolmogorov's theory of probability~\cite{Kolmogorov}.

\section{Theoretical formulation} 

\subsection{The Probability Distributions}

Both the weak and the strong LLNs are proven without imposing any assumptions on the sample space $\Omega$ and on the probability measure ${\mathbb P}$. 
Instead one only assumes the existence of the first moment, in the case of the weak law, and both the first and the second moments in the case of the strong law \cite{Varadhan}.

However, exploring the other extreme, namely the limit of the number of outcomes becoming very large ($M\rightarrow \infty$) subject to the sample size $N$ being fixed, does in fact require the knowledge of both the sample space and the probability measure. 
Hence, we have to formulate assumptions about $\Omega$ and ${\mathbb P}$.

Consider a finite discrete random variable $X$ with $(M+1)$ possible outcomes $\Omega = \lbrace m_j \rbrace_{j=0}^{M}$, with associated probabilities $\lbrace \bar{\mathbb P}_j \rbrace_{j=0}^{M}$ with $M \in \mathbb N \geq 1$.\footnote{There are $(M+1)$ possible outcomes as we index from $0$ to $M$. For compactness we will continue to use $M$ as the number of outcomes. Indeed, as $M$ becomes very large, there is little difference between $M$ and $(M+1)$.} For our purposes it is necessary to define the probability distribution as an explicit function of the outcomes, $m_j$. Thus, we use a polynomial representation given by:
\begin{equation}
\bar{{\mathbb P}}_j := \frac{{\mathbb P}_j}{\mathfrak N} = 
\frac{1}{\mathfrak N}\sum\limits_{n=0}^{M} \tilde{a}_n {m_j}^n \quad \mbox{for $j=0,\dots,M$}
\label{eq:ProbSpace_gen}
\end{equation}
Here the numbers $\left({\mathbb P}_j\right)_{j=0}^M$ are like ``unnormalized'' probabilities and the normalization factor, ${\mathfrak N}$, is given by:
\begin{equation}
{\mathfrak N} := \sum\limits_{j=0}^M {\mathbb P}_j  \label{eq:norm} 
\end{equation}
The coefficients, $\left\{\tilde{a}_n\right\}_{n=0}^M$, are real numbers which can be computed by inverting the linear relationship in (\ref{eq:ProbSpace_gen}) and thus inverting the matrix $\underline{\underline{J}}_{j, n} = {m_j}^n$. 
As $\underline{\underline{J}}$ is the VanderMonde matrix~\cite{Vandermonde} this can always be done; thus, all finite discrete probability distributions can be represented in this form.

We hope to look at random variables with more arbitrary outcomes in future work but for this paper, we look at a specific set of outcomes, $m_j = j \mbox{, }j \in \lbrace 0, 1, \ldots M \rbrace$. In this case, equation (\ref{eq:ProbSpace_gen}) becomes:
\begin{equation}
\bar{{\mathbb P}}_j := \frac{{\mathbb P}_j}{\mathfrak N} = 
\frac{1}{\mathfrak N}\sum\limits_{n=0}^{M} \tilde{a}_n j^n \quad \mbox{for $j=0,\dots,M$}
\label{eq:ProbSpace}
\end{equation}
We then use the fact that ${\mathbb P}_0 = \tilde{a}_0$, and call the quantities $\left(\frac{{\mathbb P}_j - {\mathbb P}_0}{j}\right)_{j=1}^M$ ``reduced probabilities''. We can write the other coefficients $\lbrace \tilde{a}_n \rbrace_{n=1}^M $ in terms of the reduced probabilities:
\begin{equation}
\tilde{a}_n := \sum\limits_{j=1}^M {\mathfrak A}_{j,n} \frac{{\mathbb P}_j -{\mathbb P}_0}{j}
\label{eq:InverseRel}
\end{equation}
for $n=1,\dots,M$. 
Here $\left({\mathfrak A}_{j,n}\right)_{j=1,n=1}^{M, M}$ is the inverse of
$\underline{\underline{J}} := \left(j^{n-1}\right)_{j=1,n=1}^{M,M}$.
Again, since $\underline{\underline{J}}$ is a Vandermonde matrix it always invertible. 

We note that any quantities that describe the sample (the sample mean, the variance of the sample mean, for example) depend explicitly 
on the coefficients in (\ref{eq:InverseRel}). We will see that, depending on those coefficients, we will obtain different
large-$M$ behaviour of those quantities and, in particular, of the variance of the sample mean.
Some of these coefficients may be zero. For our analysis it is useful to use $s > 0$ for the order of the polynomial in (\ref{eq:ProbSpace}) such that $\tilde{a}_s \ne 0$ and $\tilde{a}_j = 0$ for all $s<j\le M$. We will see that the variance of the sample mean scales differently with the number of outcomes depending on $s$. Both $s$ and the coefficients $\tilde{a}_n$ may vary as $M$ varies.

Recall that the variance of the sample mean depends on the variance of the underlying random variable. To look at the behaviour of the variance for different $M$, we need to find a closed form expression
for the inverse Vandermonde matrix. 
Furthermore we believe that the closed form expression for the inverse will be useful for other mathematical problems, like polynomial least square fitting, Lagrange interpolation polynomials \cite{HoffmanKunze}, and reconstruction of a statistical distribution from the moments of the distribution \cite{VonMises}.

The expression for the inverse was, to the best of our knowledge, previously unknown.
We give it below.  The inverse 
$\left({\mathfrak A}_{j,n}\right)_{n=1,j=1}^{M,M}$ reads:
\begin{eqnarray}
\lefteqn{{\mathfrak A}_{j,n} = \frac{(-1)^{j+n}}{(n-1)!(M-n)!}
\cdot
\left(\sum\limits_{p=0}^{M-j} P_{j+p}(M) (-n)^p\right)}
\label{eq:InverseVandermonde} \\
\!\!&&=\!\!
(-1)^{j+n} C^M_n
\sum\limits_{p=0}^{j-1}
\left(
\sum\limits_{1\le q_1 < \cdots < q_p \le n-1} \prod\limits_{l=1}^{p} \frac{1}{q_l}
\right)
\left(
\sum\limits_{n+1\le q_{p+1} < \cdots < q_{j-1} \le M} \prod\limits_{l=p+1}^{j-1} \frac{1}{q_l}
\right)
\label{eq:InverseVandermonde1} 
\end{eqnarray}
Here the polynomials, $P_i(j) $, satisfy the following recursion relations:
\begin{equation}
P_{M-j}(M)-P_{M-j}(M-1) = -\sum\limits_{p=1}^{j} P_{M-j+p}(M) (-M)^p
\label{eq:PolynomialsRecurs}
\end{equation}
for $j=1,\dots,M-1$ with $P_M(M) = 1$.
The lowest ten polynomials ($j=1,\dots,10$) are listed in (\ref{eq:Pol1})-(\ref{eq:Pol10}) in Appendix C.
Therein attached is also a piece of Mathematica code that tests the validity of those expressions.
The proof of (\ref{eq:InverseVandermonde}) is given in Appendix B.
Note that rows with low ($1$, $2$, $3$) and with high ($M$, $M-1$, $M-2$) indices  have a particularly 
simple form; the complexity of the expression rises when the row index tends towards $M/2$.
Setting $j=1,2,3$ in (\ref{eq:InverseVandermonde1}) and $j=M,M-1,M-2,M-3$ in (\ref{eq:InverseVandermonde}) we obtain:
\begin{eqnarray}
\lefteqn{{\mathfrak A}_{1,n} = (-1)^{1+n} C^M_n}  \\
&&\!\!\!\!
{\mathfrak A}_{2,n} = 
(-1)^1 {\mathfrak A}_{1,n} \left(\sum\limits_{q=1}^{n-1} \frac{1}{q} (1_{1\le q \le n-1} + 1_{n+1\le q \le M})\right) \\
&&\!\!\!\!
{\mathfrak A}_{3,n} = (-1)^2 {\mathfrak A}_{1,n} 
\left( \sum\limits \frac{1}{q_1 q_2} 1_{n+1\le q_1<q_2\le M} \right. \nonumber \\
&&
\left. + \sum\limits \frac{1}{q_1 q_2} 1_{1\le q_1 \le n-1} 1_{n+1 \le q_2 \le M} 1_{1\le q_1 < q_2 \le n-1}\right) \\
&&\vdots \nonumber \\
&&\!\!\!\!
\frac{{\mathfrak A}_{M-3,n}}{{\mathfrak A}_{M,n}(-1)^{-3}} =   
  \left( \frac{(M-2)(M-1)M^2(M+1)^2}{48}\right. \nonumber \\ 
&&
 \left. -\frac{n}{24}(M-1)M(M+1)(2+3 M) + \frac{n^2}{2} M(M+1) -n^3\right)\\
&&\!\!\!\!
\frac{{\mathfrak A}_{M-2,n}}{{\mathfrak A}_{M,n}(-1)^{-2}} =  
\left(\frac{(M-1)M(M+1)(2+3 M)}{24} - \frac{n}{2} M(M+1)  + n^2\right) \\
&&\!\!\!\! 
\frac{{\mathfrak A}_{M-1,n}}{{\mathfrak A}_{M,n}(-1)^{-1}} =   
  \left(\frac{M(M+1)}{2}  - n  \right) \\
&&\!\!\!\!
{\mathfrak A}_{M,n} = \frac{(-1)^{M+n}}{(n-1)!(M-n)!}
\end{eqnarray}

\noindent{\bf Note:} The inversion of the Vandermonde matrix has applications in 
control theory \cite{Tou,deBrule,Reis}, in signal processing problems \cite{Neagoe}
and in systems theory \cite{Reis,Kan,Goeknar,Wertz}.
In order to solve the problem one typically makes use of the Lagrange interpolation formula 
and finds the rows of the inverse matrix 
by computing the coefficients of the Lagrange interpolation polynomials related to the matrix elements \cite{Tou}.
An alternative approach was proposed in \cite{Kaufman} where the inverse matrix is expressed 
as a product of two matrices
one of which is diagonal and the elements of the other one are given through certain recursion relations. 
In this way the total number of operations needed to invert is reduced from $O(M^3)$ to $O(M^2)$.
Finally in \cite{Neagoe} one expresses the elements of the inverse through totally symmetric polynomials
of the elements of the original matrix.

Our method of computing the inverse (presented in Appendix B) is a way of actually 
solving the recursion relations from \cite{Kaufman} or summing up
the totally symmetric polynomials analytically if the elements of the original matrix are certain real or complex 
powers of a constant or of an arithmetic progression.
Indeed, even though the elements of the original matrix are first powers of an arithmetic progression in our work i.e.\ $x_j=j$,  the manipulations
(\ref{eq:MatrixInverse})-(\ref{eq:MatrixInverse1a}) can also be done analytically in the generic case of arbitrary powers, with little effort.
Thus, we believe, that our method is superior to the methods known in the literature
and can be applied to produce more efficient numerical algorithms of use in the areas described above.

Now we fix $N \ge 1$, we draw a sample of size $N$ from the population described in (\ref{eq:ProbSpace}) and we conjecture that the quantities that describe the sample (estimators of population parameters) diverge from those that describe the population, if the number of outcomes $M$ becomes very large.
To be specific, we analyze the variance of the sample mean, and we show that, except for very ``unusual distributions'', it increases monotonically with $M$, if $M$ is big enough, which implies that, for fixed $N$, the sample mean diverges from the population mean when $M \to \infty$.
In most computations below we will use asymptotic limits ($M\to \infty$) rather than exact results.
However, deriving exact results is not an essential difficulty; it amounts to performing more work, which is not necessary for our purposes.

\subsection{The variance of the sample mean} 

The variance of the sample mean reads:
\begin{equation}
\sigma^2_{\overline{X}_N} := \mbox{var}\left[ \frac{\sum\limits_{i=1}^N X_i}{N} \right] = \frac{\sigma^2}{N} := \frac{\left<\left(X - \left<X\right>\right)^2\right>}{N} = \frac{\left<X^2\right> - \left<X\right>^2}{N}
\label{eq:SampleVardef}
\end{equation}
It depends on the variance, $\sigma^2$, of the random variable. Using equation (\ref{eq:ProbSpace}) and Faulhaber's formula (\ref{eq:Faulhaber}), this variance can be written as:
\begin{eqnarray}
\sigma^2 &=& 
\frac{1}{{\mathfrak N}^2}
\left[
{\mathfrak N} \left(\sum\limits_{j=0}^M j^2 {\mathbb P}_j\right)
-
\left(\sum\limits_{j=0}^M j {\mathbb P}_j\right)^2
\right]
=
\frac{\sum\limits_{n=0}^{2{M}+4} \beta^{(1)}_n
(M+1)^n
}{
\sum\limits_{n=0}^{2{M}+2} \beta^{(2)}_n
(M+1)^n
}
\label{eq:MyVariance}
\end{eqnarray}
where
\begin{eqnarray}
\beta^{(1)}_n &:=& 
\sum\limits_{n=n_1+n_2+4-(k_1+k_2)} 1_{n_1\le M}1_{n_2\le M}
\tilde{a}_{n_1} \tilde{a}_{n_2} B_{k_1} B_{k_2} \nonumber\\ 
&&
\left(
\frac{C^{n_2+3}_{k_2}  1_{k_2\le n_2+3}}{n_2+3} \frac{C^{n_1+1}_{k_1}1_{k_1\le n_1+1}}{n_1+1}  
-
\frac{C^{n_2+2}_{k_2}  1_{k_2\le n_2+2}}{n_2+2} \frac{C^{n_1+2}_{k_1}1_{k_1\le n_1+2}}{n_1+2}  
\right)\label{eq:CoeffsRatioFct0}\\
\beta^{(2)}_n &:=& 
\sum\limits_{n=n_1+n_2+2-(k_1+k_2)} 1_{n_1\le M}1_{n_2\le M}
\tilde{a}_{n_1} \tilde{a}_{n_2} B_{k_1} B_{k_2} \nonumber \\
&& \times \left( \frac{C^{n_1+1}_{k_1}1_{k_1\le n_1+1}}{n_1+1}\right) \left( \frac{C^{n_2+1}_{k_2} 1_{k_2\le n_2+1}}{n_2+1} \right) 
\label{eq:CoeffsRatioFct}
\end{eqnarray}
where $n_1, n_2, k_1, k_2 \in \mathbb{N}$. The numbers, $B_{k_i}$, are the Bernoulli numbers.

Thus both the variance $\sigma^2$ and the variance of the sample mean are rational functions of the number of outcomes $M$.
Since we are only interested in the large-$M$ behaviour of the variance of the sample mean, we do not need to simplify expressions (\ref{eq:CoeffsRatioFct0}) and (\ref{eq:CoeffsRatioFct}) for the coefficients but instead we only need to work out leading order terms in both the numerator and the denominator in (\ref{eq:MyVariance}).
Thus, we assume that the order of the polynomial in (\ref{eq:ProbSpace}) is $s>0$, meaning
that the coefficent $\tilde{a}_s$ is non-zero and the last (M-s) coefficients, $\lbrace \tilde{a}_n \rbrace_{n=s+1}^{M}$ in the expansion (\ref{eq:ProbSpace}), are zero.
Then the coefficients for the highest order terms in the numerator read:
\begin{eqnarray}
\beta^{(1)}_{2 s+4} &=&  \frac{(\tilde{a}_s)^2}{(s+1)(s+3)(s+2)^2} \\
\beta^{(1)}_{2 s+3} &=&  \frac{2 \tilde{a}_{s-1} \tilde{a}_s}{s(s+1)(s+2)(s+3)} -  \frac{(\tilde{a}_s)^2}{(s+1)(s+2)(s+3)} \\
\beta^{(1)}_{2 s+2} &=&   \frac{-\tilde{a}_s \tilde{a}_{s-2}}{s(s+1)^2 (s+2)} +   
\frac{\tilde{a}_{s-1} \tilde{a}_{s-1}}{s(s+1)^2 (s+2)} +  \frac{3 \tilde{a}_s \tilde{a}_{s-2}}{(s-1)s(s+2)(s+3)}
\\\nonumber &&
+  \frac{\tilde{a}_s \tilde{a}_{s-1} 2(3+2s)}{s(s+1)(s+2)(s+3)} 
+  \frac{(\tilde{a}_s)^2 (2s+3)}{(s+1)(s+2)(s+3)} \\
\beta^{(2)}_{2 s+2} &=&  \frac{(\tilde{a}_s)^2}{(s+1)^2} \\
\beta^{(2)}_{2 s+1} &=&  \frac{2\tilde{a}_{s-1} \tilde{a}_s}{s(s+1)} -  \frac{(\tilde{a}_s)^2}{(s+1)} 
\end{eqnarray} 

Dividing the polynomials in (\ref{eq:CoeffsRatioFct0}) and (\ref{eq:CoeffsRatioFct}) by one another we obtain the following expression for $\sigma^2$, for large $M$:
\begin{eqnarray}
 \sigma^2 =
(M+1)^2 \frac{(s+1)}{(s+3)(s+2)^2}
 +
(M+1) \left(2 \frac{\tilde{a}_{s-1}}{\tilde{a}_s} - s\right) \frac{(s+1)}{s (s+2)^2(s+3)}
+O(1)
\label{eq:Variance}
\end{eqnarray}
Subbing this into equation (\ref{eq:SampleVardef}) will also give an expression for the variance of the sample mean as a function of $M$.

\subsection{The Large $M$ behaviour of the Variance}\label{ss:Large_M}

We now wish to examine whether $\sigma^2$ converges or diverges for large $M$. Equation (\ref{eq:Variance}) can be rewritten as:
\begin{equation}
\sigma^2 = 
\frac{(s+1)}{(s+3)(s+2)^2} \left[M^2 + 2M \left(\frac{\tilde{a}_{s-1}}{\tilde{a}_s}\right) \left(\frac{1}{s}\right) + M + 2\left(\frac{\tilde{a}_{s-1}}{\tilde{a}_s}\right) \left(\frac{1}{s}\right) \right]
\label{eq:Variance1}
\end{equation}
As noted previously, it is possible for both the order, $s$, and the coefficients, $\tilde{a}_n$, to vary with $M$. Looking at the term in the square brackets, we notice that the first term ($M^2$) will always dominate the third term ($M$) and the second term ($2M \left(\frac{\tilde{a}_{s-1}}{\tilde{a}_s}\right) \left(\frac{1}{s}\right)$) will always dominate the fourth term ($2\left(\frac{\tilde{a}_{s-1}}{\tilde{a}_s}\right) \left(\frac{1}{s}\right)$). Therefore we want to consider the large $M$ behaviour due to the first and second terms. The exact behaviour of the second term will depend on the distribution concerned, but we can identify four broad cases depending on the behaviour of the ratio $\left(\frac{\tilde{a}_{s-1}}{\tilde{a}_s}\right)\left(\frac{1}{s}\right)$ in the large $M$ limit:
\begin{enumerate}
\item
$\displaystyle\lim_{M \to \infty} \left|\frac{\tilde{a}_{s-1}}{\tilde{a}_s}\right| \left(\frac{1}{s}\right) < \infty$
\item
$\displaystyle\lim_{M \to \infty} \left|\frac{\tilde{a}_{s-1}}{\tilde{a}_s}\right| \left(\frac{1}{s}\right) = \infty$ and  $2M\left| \frac{\tilde{a}_{s-1}}{\tilde{a}_s}\right| \left(\frac{1}{s}\right) \ll M^2$ 
\item
$\displaystyle\lim_{M \to \infty} \left(\frac{\tilde{a}_{s-1}}{\tilde{a}_s}\right) \left(\frac{1}{s}\right) = \infty$ and $2M\left( \frac{\tilde{a}_{s-1}}{\tilde{a}_s}\right) \left(\frac{1}{s}\right) \gtrsim  M^2$ 
\item
$\displaystyle\lim_{M \to \infty}\left(\frac{\tilde{a}_{s-1}}{\tilde{a}_s}\right) \left(\frac{1}{s}\right) = -\infty$ and $2M\left( \frac{\tilde{a}_{s-1}}{\tilde{a}_s}\right) \left(\frac{1}{s}\right) \lesssim -M^2$ 
\end{enumerate}
The first case corresponds to when the limit of the ratio is finite. The second case is when the ratio grows in absolute value with $M$ but at a rate slower than $M$. The third case is when the ratio grows with $M$ at a rate faster than $M$. Finally, the fourth case is when the ratio grows in absolute value with $M$ but is negative. We now consider each case separately.

\subsubsection{Case 1} 

In this case, we do not need to worry about the coefficients, and how they depend on $M$ to understand the large $M$ behaviour of the variance. It is clear that only the first term from equation (\ref{eq:Variance1}) matters as $M \to \infty$. Thus we have:
\begin{equation}
\sigma^2 \approx \frac{(s+1)}{(s+3)(s+2)^2} M^2 \label{eq:case1}
\end{equation}

Now, let us consider three examples. Firstly, $s$ does not depend on $M$; secondly, 
$s \simeq \sqrt{M}$, i.e.\ $\lim\limits_{M\rightarrow \infty} s/\sqrt{M} = \gamma$ where $0 < \gamma < \infty$; and finally,  $s \simeq M$, i.e.\ $\lim\limits_{M\rightarrow \infty} s/M = \tilde{\gamma}$ where $0 < \tilde{\gamma} < \infty$.
From (\ref{eq:case1}) we see that in all three instances, both the variance and the variance of the sample mean (for given $N$)
diverge as the second, the first and the zeroth power of the number of outcomes respectively. That is, we have: 
\begin{eqnarray}
\lim_{M\rightarrow \infty} \frac{\sigma^2}{M^2} &=&  \frac{(s+1)}{(s+3)(s+2)^2} < \infty \label{eq:AccDiv}\\
\lim_{M\rightarrow \infty} \frac{\sigma^2}{M^1} &=&  \frac{1}{\gamma^2} \label{eq:Div}\\
\lim_{M\rightarrow \infty} \frac{\sigma^2}{M^0} &=&  \frac{1}{\gamma_1^2}\label{eq:DecDiv}
\end{eqnarray}

We now go on to provide three explicit examples of such distributions. In all three examples, for fixed sample size $N$, the sample mean, or the estimator of the population mean, diverges from the population mean, yet at different rates.
We describe these examples as \textit{\textbf{Accelerating Divergence}}, \textit{\textbf{Divergence}}, and \textit{\textbf{Decelerating Divergence}}, respectively.

\begin{figure}[hbt]
\hbox{
\psfig{figure=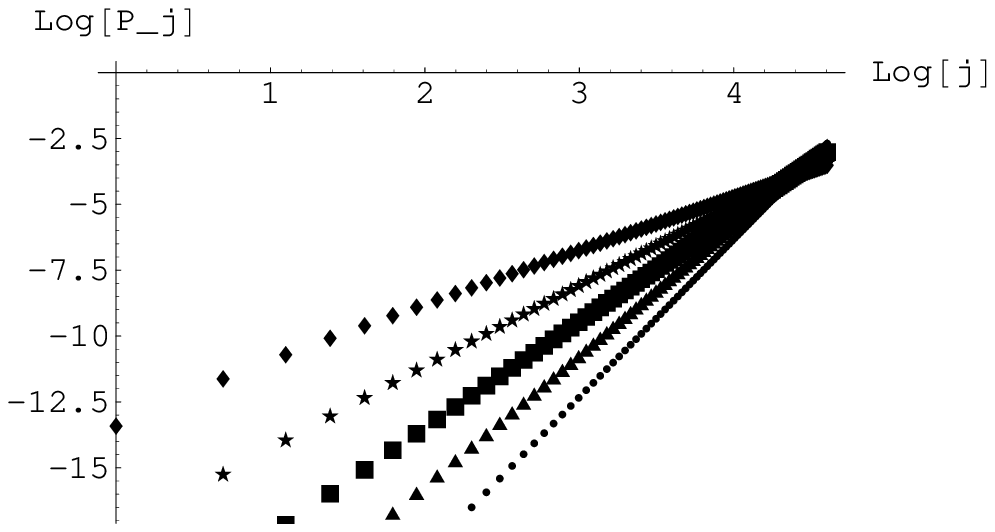,width=0.49\textwidth}
\psfig{figure=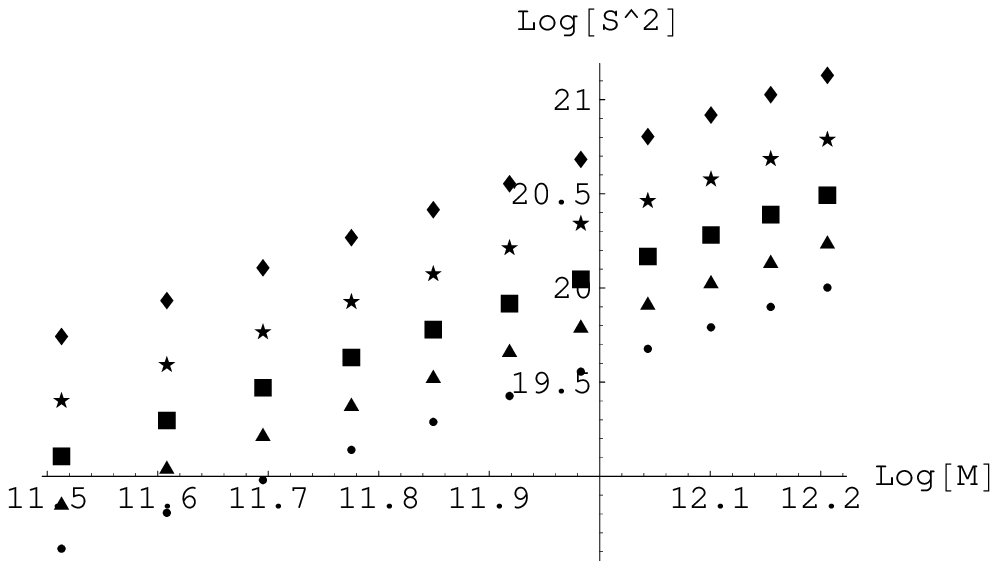,width=0.49\textwidth}
}
\caption{The probability distributions in equation (\ref{eq:ProbSpace}) for different 
values of the order of the polynomial $s=2,\dots,6$ (left)
and the variance of that distribution as a function of the number of outcomes $M=100000,\dots,200000$ (right).
A second order polynomial was fitted by least square regression to each of the data sets on the right. (This was done to confirm the relation in equation (\ref{eq:case1})).
The parameters of that fit are
(-3.2831 (-3.2834), 1.998, 0.000087), (-3.6240 (-3.6243), 1.998, 0.000097), (-3.9195 (-3.9199), 1.998, 0.000138), (-4.17872 (-4.1795), 1.9959, 0.00023), (-4.4096 (-4.4102), 1.99677, 0.00018) for $s=2,\dots,6$ respectively.
The theoretical values of the first parameter (the intercept) are given in brackets and read
$\log((s+1)/((s+3)(s+2)^2)$.
The theoretical values of the second and the third parameter are two and zero respectively. 
\label{fig:ProbVar}
}
\end{figure}

\paragraph{Accelerating divergence:}
We generated a sequence of non-zero real parameters $\left(\tilde{a}_n\right)_{n=1}^s$ and computed the normalized probabilities 
$\left(\bar{\mathbb P}_j\right)_{j=0}^M$ from (\ref{eq:ProbSpace}).
We run over $s=2,\dots,6$ and for each value of $s$ we plot both the normalized probabilities, for $M=100$, as a function of $j$, and the variance of the distribution as a function of the number of outcomes $M$, in the double logarithmic scale (see Figure \ref{fig:ProbVar}).
It is readily seen that the variance behaves asymptotically as the second power of the number of outcomes, and, in addition the asymptotic behaviour is attained quite quickly.

\noindent{\bf Note:}
We reiterate that this type of divergence is obtained for every distribution whose unnormalized probability distribution can be represented by a polynomial of order $s$ that does not depend on $M$, when $M$ is large.

\begin{figure}[hbt]
\hbox{
\psfig{figure=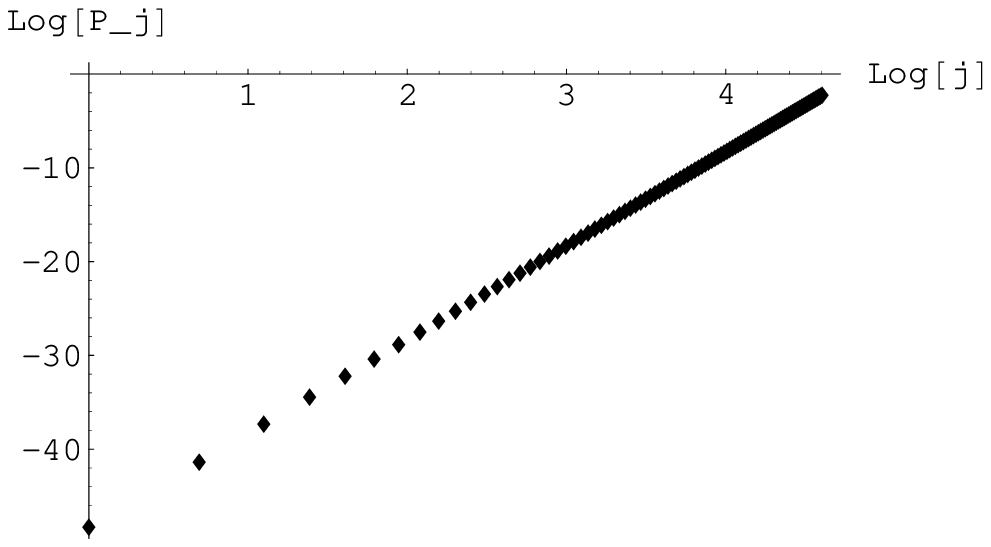,width=0.49\textwidth}
\psfig{figure=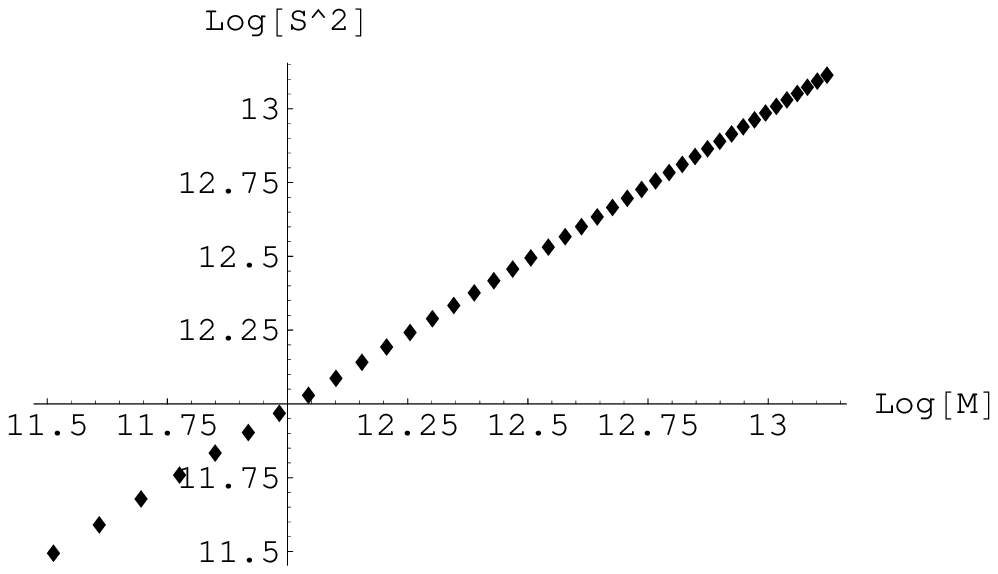,width=0.49\textwidth}}
\caption{The probability distribution ${\mathbb P}_j = j^{\sqrt{M}}$ (left) and the variance of that distribution as a function of the number of outcomes $M=100000,500000$ (right)
Here, the parameters of the linear regression second order polynomial fit are $(-0.3\pm 0.0027,1.04\pm 0.0004,0.002\pm 0.000017)$ compared to the theoretical parameters $(0.0,1.0,0.0)$,
\label{fig:ProbVar1}
}
\end{figure}

\paragraph{Divergence:} 
Here we repeat the procedure from the previous point with one difference, namely that the number of parameters depends on the number of outcomes like equation (\ref{eq:Div}).
For the sake of simplicity let us take $\tilde{a}_n = \delta_{n,\gamma \sqrt{M}}$ where $\gamma>0$.
Then, asymptotically, the unnormalized probabilities ${\mathbb P}_j$, and the norm ${\mathfrak N}$, read ${\mathbb P}_j = j^{\gamma \sqrt{M}}$, ${\mathfrak N} = \frac{M^{\gamma \sqrt{M} + 1}}{(\gamma\sqrt{M}+1)}$ respectively.
The probability distribution is plotted in the first graph of Figure~\ref{fig:ProbVar1}.
The first and second moments, and the variance read:
\begin{eqnarray}
\left<j\right> &=& \sum\limits_{j=0}^M j j^{\gamma \sqrt{M}} \left(\frac{\gamma \sqrt{M} + 1}{M^{\gamma \sqrt{M}+1}}\right)
=
M \left( \frac{\gamma \sqrt{M}+1}{\gamma \sqrt{M}+2} \right) \label{eq:FirstMoment}\\
\left<j^2\right> &=& M^2 \left( \frac{\gamma \sqrt{M}+1}{\gamma \sqrt{M}+3} \right) \label{eq:2ndMoment}\\
\sigma^2 &=& \left<j^2\right> - \left<j\right>^2 = M^2 \left( \frac{\gamma \sqrt{M}+1}{(\gamma \sqrt{M}+3)(\gamma \sqrt{M}+2)^2} \right)
\mathop{=}_{M\rightarrow \infty} \frac{1}{\gamma^2} M
\end{eqnarray}
In Figure \ref{fig:ProbVar1} we plot the variance of the distribution as a function of the number of outcomes $M$.
Now the variance diverges asymptotically as the first power of the number of outcomes. Here, however,
the asymptotic behaviour is attained much slower than in the previous case.

\noindent{\bf Note:} We emphasize that this type of divergence is characteristic for every distribution whose probability function can be represented by a polynomial of order $s$  that behaves like the square root of the number of outcomes $M$, for large $M$.
Indeed, setting $s=s_0\sqrt{M}$ in (\ref{eq:Variance1}) we get:
\begin{equation}
\sigma^2 \mathop{=}_{M\rightarrow\infty} M \frac{1}{{s_0}^2} + O(1)
\end{equation}

\begin{figure}[thb]
\hbox{
\psfig{figure=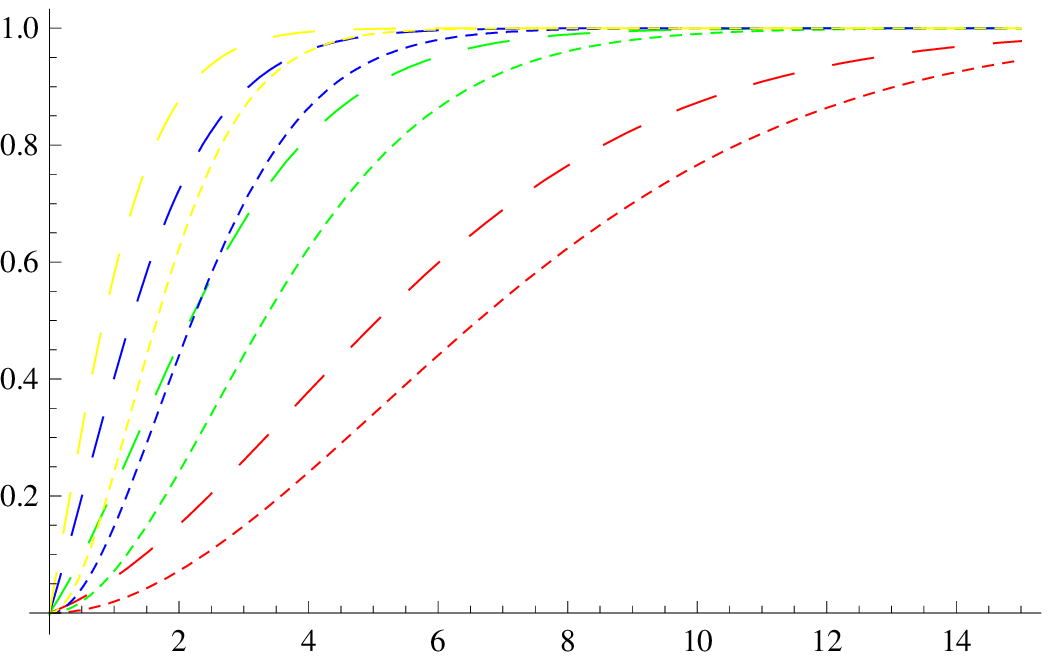,width=0.48\textwidth}
\psfig{figure=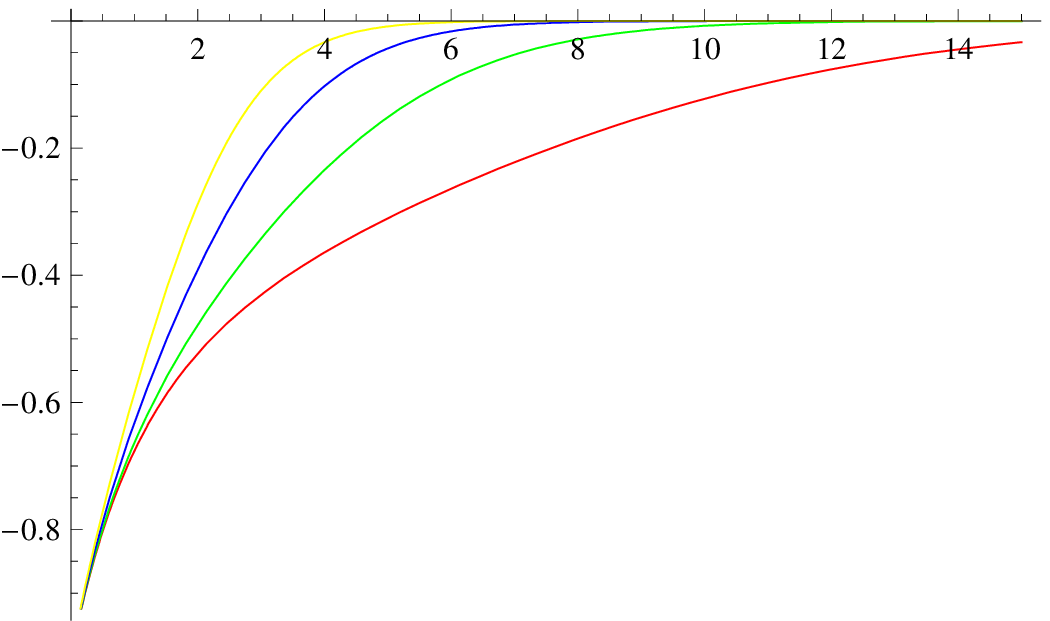,width=0.48\textwidth}
}
\caption{On the left: The approximated variance (short dash) versus the exact variance (long dash) of an exponential distribution as a function of the number of outcomes for different values of $\alpha=1.0,1.5,2.0$ (from top to bottom).
On the right: The difference of the approximated and the exact variance  divided by the exact variance as a function of the number of outcomes.  
\label{fig:SampleVariance}}
\end{figure}

\paragraph{Decelerating Divergence:}
Finally, we consider an exponential distribution 
\begin{equation}
\bar{{\mathbb P}}_j = \exp(-j\alpha) \left( \frac{1- \exp(-\alpha)}{1 - \exp(-(M+1)\alpha)}\right)
\label{eq:ExpDistr}
\end{equation}
for $j=0,\ldots,M$ and $\alpha > 0$. 
For large $M$,
the distribution can be well approximated by its Taylor expansion truncated at the $M$\textsuperscript{th} order
and as such, it would correspond to the polynomial form in equation (\ref{eq:ProbSpace}).
Then $\tilde{a}_n = (-\alpha )^n/n!$  
and, in the large $M$ limit, the variance takes the following form:
\begin{equation}
\sigma^2 = \frac{1}{\alpha^2} 
\left( \frac{1 - \exp(-\alpha M)(1+\alpha M+\frac{(\alpha M)^2}{2}+\frac{(\alpha M)^3}{6})}{1 - \exp(-\alpha M)(1+\alpha M)}\right)
\mathop{=}_{M\rightarrow \infty}
\frac{1}{\alpha^2} 
\label{eq:SampleVarExp}
\end{equation}
Recall that equation (\ref{eq:SampleVarExp}) is only an approximation. In order to check the validity of that
approximation, we plot the approximated variance of the sample mean $\sigma^2_{\overline{X}_N}=\sigma^2/N$ along with the exact result\footnote{We have derived the exact variance of the sample mean using MATHEMATICA. The expression is lengthy and cumbersome and, in our opinion, it does not bring much to quote it here.} in Figure \ref{fig:SampleVariance}. 
We can see from the figure that the results do not differ by more than $1\%$ for $M>10$.

Now the variance behaves asymptotically as the zeroth power of the number of outcomes and this behaviour is attained much faster than in the previous case.

\noindent{\bf Note 1:} Again we stress that this type of divergence is characteristic for every distribution whose probability 
function is a polynomial whose order $s$ is proportional to the number of outcomes $M$.
Indeed, setting $s=s_ M$ in (\ref{eq:Variance}) we get:
\begin{equation}
\sigma^2 \mathop{=}_{M\rightarrow\infty}  \frac{1}{{s_0}^2} + O(\frac{1}{M})
\end{equation}

From the above discussion it is clear that the condition:
\begin{equation}
\lim_{M \to \infty} \left|\frac{\tilde{a}_{s-1}}{\tilde{a}_s}\right| \left(\frac{1}{s}\right) < \infty \label{eq:Case1_Constraint}
\end{equation}
is sufficient for the variance to exhibit divergent behaviour in the large $M$ limit.
It is however an open question to specify both \textit{\textbf{necessary and sufficient conditions}} for distributions to exhibit the behaviour specified above.
One way of doing that, is to require that the last $(M-s)$ coefficients $\left(\tilde{a}_n\right)_{n=s+1}^M$ are equal to zero.
This, from (\ref{eq:InverseRel}), yields $(M-s)$ linearly independent linear equations for the unknown probabilities, given by: 
\begin{equation}
\tilde{a}_n = \sum\limits_{j=1}^M {\mathfrak A}_{j,n} \frac{{\mathbb P}_j - {\mathbb P}_0}{j} = 0
\label{eq:SysEq}
\end{equation}
for $n=s+1, s+2, \ldots, M$.
The solution to (\ref{eq:SysEq}) is given in (\ref{eq:ReducedProbabilitiesFinal1}) 
along with (\ref{eq:ModesProbs})
and its proof is in Appendix D.

\paragraph{Summary of Case 1} An $s$-parameter family of unnormalized distributions 
is given by:
\begin{equation}
{\mathbb P}_j = {\mathbb P}_0 + \sum\limits_{q=1}^s {\mathfrak V}_{q,j} \frac{j}{q} ({\mathbb P}_q - {\mathbb P}_0)
\label{eq:ReducedProbabilitiesFinal11}
\end{equation}
for $j=s+1,\dots,M$ with the quantities ${\mathfrak V}_{q,j}$ being defined in (\ref{eq:ModesProbs}). 
A probability distribution satisfying (\ref{eq:ReducedProbabilitiesFinal11}) and the constraint (\ref{eq:Case1_Constraint}) has a property that its variance $\sigma^2$, along with the variance of the sample mean $\sigma^2_{\overline{X}_N}$ (for a given $N$), exhibits an asymptotic behaviour as in (\ref{eq:Variance1}).
This implies that the large-$M$ behaviour of the variance of the sample mean can range in a continuous fashion from a quadratic divergence $(\sigma^2_{\overline{X}_N} \simeq M^2)$ to asymptotic convergence $(\sigma^2_{\overline{X}_N} \simeq M^0)$, for fixed sample size, $N$.
Furthermore, if we assume the order of the polynomial in (\ref{eq:ProbSpace}), $s$, has a power law dependence on $M$ for large $M$, of the form  $s \simeq M^{\alpha}$, with $0 \le \alpha \le 1$, then, from (\ref{eq:Variance1}), and for fixed $N$ the variance of the sample mean behaves as $\sigma^2_{\overline{X}_N} \simeq M^{2-2\alpha}$.

\begin{figure}[tbh]
\hbox{
\psfig{figure=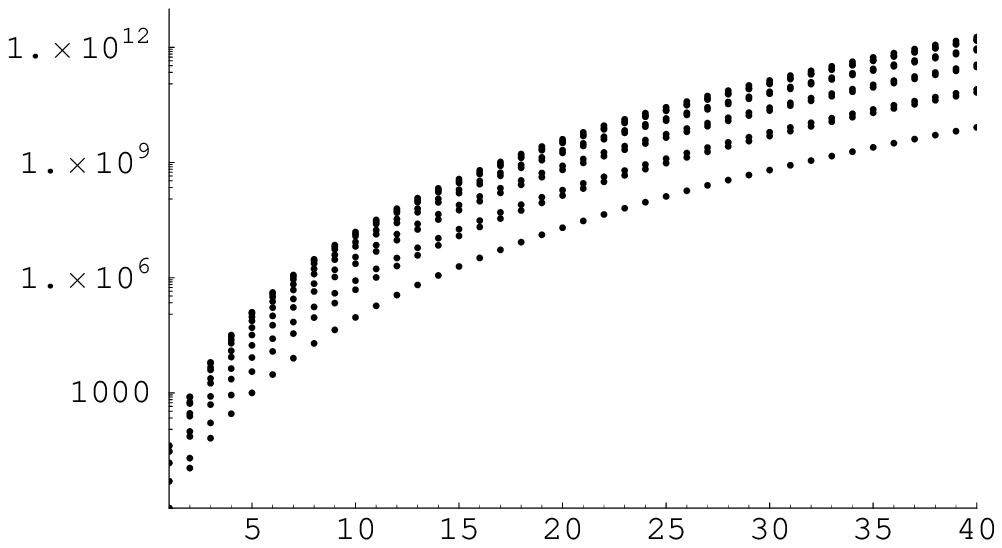,width=0.45\textwidth}
\psfig{figure=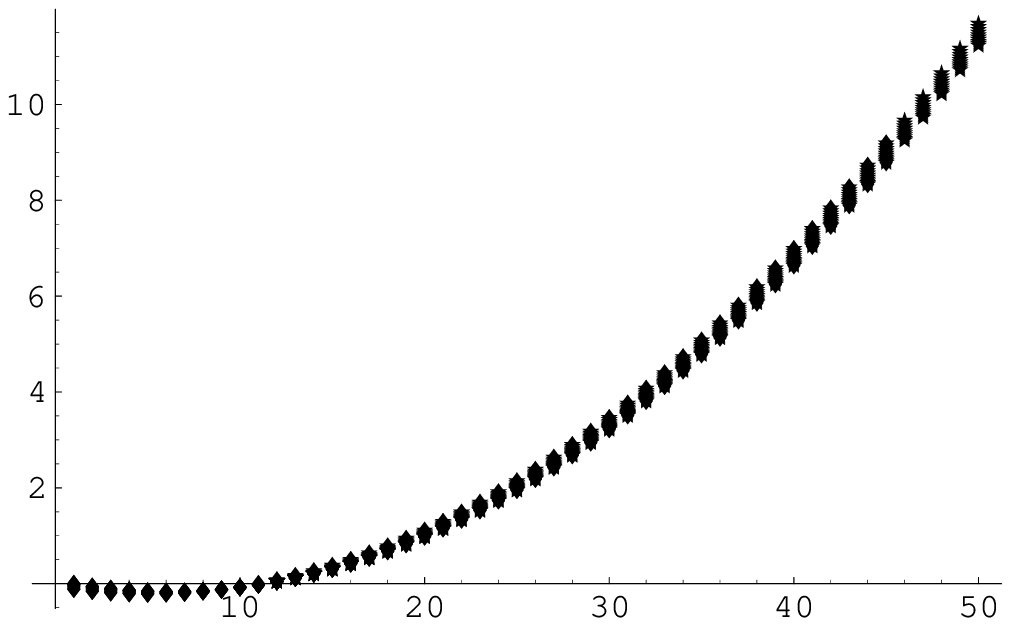,width=0.45\textwidth}
}
\caption{On the left: All probability distributions for case 1. \label{fig:ProbDistrsModes}
Here $M=50$ and $s=10$. We plot the quantities $\left({\mathfrak V}_{q,j} \frac{j}{q}\right)_{j=s+1}^M$ for $q=1,\dots,s$ and the distributions are constructed as in (\ref{eq:ReducedProbabilitiesFinal11}).
On the right: The dependence of the variance of the probability distributions shown in the left of Figure \ref{fig:ProbDistrsModes} on the number of outcomes. 
Here $s=10$ is fixed and $M=s+2,\dots,50$. 
The graphs are overlayed with a polynomial of second order fitted by least square regression.
\label{fig:VarianceNumPlots}.
The leading order coefficients read $0.00603$, $0.00606$, $0.00607$, $0.00609$,  $0.00609$, $0.00608$, $0.00607$, $0.00601$, $0.00593$, and $0.00587$ compared to the theoretical value $0.00587$ (compare equation (\ref{eq:Variance})).}
\end{figure}

We did some numerical testing of our results. For particular values of $M=50$ and $s=10$ we have found the solutions (\ref{eq:ModesProbs}) by solving numerically  equations
(\ref{eq:SysEq}). 
Subsequently we found the unnormalized probabilities from (\ref{eq:ReducedProbabilitiesFinal11}) by taking arbitrarily ${\mathbb P}_0=0$ and ${\mathbb P}_q=(-1)^q$ for $q=1,\dots,s$.
We plot the results in Figure \ref{fig:ProbDistrsModes}.

In addition, for a particular value of $s$ ($s=10$), and for $M=12,\dots,50$, we computed numerically the variance of all the $s$ different probability distributions in (\ref{eq:ReducedProbabilitiesFinal11}) and we plotted the results, as a function of $M$ in Figure \ref{fig:VarianceNumPlots} as well.
As we can see, the variance displays a quadratic dependence on the number of outcomes and the leading order coefficient fits in well with the result in (\ref{eq:Variance1}).

\subsubsection{Case 2}

Recalling section \ref{ss:Large_M}, case 2 includes all probability distributions that satisfy the following conditions:
\begin{eqnarray}
\lim_{M \to \infty} \left|\frac{\tilde{a}_{s-1}}{\tilde{a}_s}\right| \left(\frac{1}{s}\right) = \infty \label{eq:Case2_Constraint1} \\
2M\left| \frac{\tilde{a}_{s-1}}{\tilde{a}_s}\right| \left(\frac{1}{s}\right) \ll M^2 \label{eq:Case2_Constraint2}
\end{eqnarray}
For this case, looking at equation (\ref{eq:Variance1}) it is clear that only the first term is important for large $M$. Thus one can use equation (\ref{eq:case1}) in the previous section and follow the discussion of divergence there. Thus case 2 distributions will also exhibit divergence in the variance as $M \to \infty$.  This implies, as before, that, for fixed $N$, the large-$M$ behaviour of the variance of the sample mean can range in a continuous fashion from a quadratic divergence $(\sigma^2_{\overline{X}_N} \simeq M^2)$ to asymptotic convergence $(\sigma^2_{\overline{X}_N} \simeq M^0)$. 

\subsubsection{Case 3}

Case 3 distributions are those that satisfy:
\begin{eqnarray}
\lim_{M \to \infty} \left(\frac{\tilde{a}_{s-1}}{\tilde{a}_s}\right) \left(\frac{1}{s}\right) = \infty \label{eq:Case3_Constraint1}\\
2M\left( \frac{\tilde{a}_{s-1}}{\tilde{a}_s}\right) \left(\frac{1}{s}\right) \gtrsim  M^2 \label{eq:Case3_Constraint2}
\end{eqnarray}
Unlike in previous cases, it is clear that now the second term in equation (\ref{eq:Variance1}) will contribute to the variance in the large $M$ limit. However, the variance and the variance of the sample mean (for fixed $N$) will still diverge as $M \to \infty$. But if the inequality in (\ref{eq:Case3_Constraint2}) is strict, then the second term in equation (\ref{eq:Variance1}) now controls the exact behaviour of the divergence and it will be faster divergence than the maximum possible ($M^2$) in the previous two cases. However, unlike in the previous cases, we do not analyze this behaviour any further.

\subsubsection{Case 4} 

In the previous three cases, for a given $N$, the variance of the sample mean increases with the number of outcomes, which might imply that this holds for every distribution. 
However this is not the case.
One can construct distributions such that the variance of the sample mean decreases with the number of outcomes.

For example, take the (unnormalized) probability distriubtion given by:
\begin{equation}
{\mathbb P}_j  = j^{M^2}\label{eq:Counter_Ex}
\end{equation}
Although this suggests a polynomial of order $M^2$, for each $M$, it can always be decomposed in a unique way into a polynomial of the form in (\ref{eq:ProbSpace}).
In this case the leading order term in the decomposition corresponds to $s=M$. 
Asymptotically, we get ${\mathfrak N} = M^{M^2+1}/(M^2+1)$ and the variance reads:
\begin{equation}
\sigma^2 = M^2 \left( \frac{1+M^2}{(3+M^2)(2+M^2)^2} \right)
\label{eq:DecayingDiverg}
\end{equation}
Hence the variance is inversely proportional to the square of the number of outcomes $M$, for large $M$. Thus we have convergent behaviour for the variance.

If we now solve for the ratio $\frac{\tilde{a}_{s-1}}{\tilde{a}_s}$  in equation (\ref{eq:Variance1}) by requiring that the variance behaves like $M^{-2}$ for large $M$, then we get:
\begin{equation}
\frac{\tilde{a}_{s-1}}{\tilde{a}_s} = \frac{1}{2}\left[ \frac{\left(-M^5 - 2M^4 + 7M^2 + 16M + 12 \right)}{\left( M(M^2 + 2M + 1 \right)} \right]\label{eq:ratio}
\end{equation}

This equation combined with the fact that $s = M$ confirms that the probability distribution in (\ref{eq:Counter_Ex}) fits under case 4 as its satisfies the conditions given by:
\begin{eqnarray}
\lim_{M \to \infty}\left(\frac{\tilde{a}_{s-1}}{\tilde{a}_s}\right) \left(\frac{1}{s}\right) = -\infty \label{eq:Case4_Constraint1} \\
2M\left( \frac{\tilde{a}_{s-1}}{\tilde{a}_s}\right) \left(\frac{1}{s}\right) \lesssim -M^2 \label{eq:Case4_Constraint2}
\end{eqnarray}

It is also possible to show that all distributions of the form ${\mathbb P}_j  = j^{M^x}$ where $x >1$, like that in equation (\ref{eq:Counter_Ex}), will have a variance that converges as $M \to \infty$. These distributions are, however, quite unusual since most probabilities are assigned to very few outcomes. It is still an open question as to how many distributions will have variances that will exhibit convergent behaviour for a large number of outcomes.

\subsection{The Percentiles}

So far we used the variance as a measure of the ``width'' of the distribution.
This may be unsatisfactory because for certain distributions the variance may not carry the relevant information - an example is the case where the variance does not exist. It is thus much more useful to compute percentiles of the distribution of the difference between the sample and the population mean, $Z := (\sum\limits_{i=1}^N X_j)/N - \mu$, and investigate the conditions under which they diverge. 
Here $\mu := \left<X\right>$ is the mean of the random variable $X$.
Recall that the $p$-percentile $z_p$ of random variable $Z$ is defined as such a value of argument such that the Cumulative Distribution Function (CDF) equals $p$\%, i.e.\ 
\begin{equation}
P\left( Z > z_p \right) = \frac{p}{100}
\end{equation}
Here $0< p < 100$.
We denote by ${\mathfrak N} := \sum\limits_{j=0}^M P_j$ the normalization constant and we compute the CDF as follows:
\begin{eqnarray}
P\left(Z > z\right) &=& \sum\limits_{j_1,\dots,j_N=0}^M 1_{\sum\limits_{p=1}^N j_p > (z+\mu) N} \prod\limits_{p=1}^N P_{j_p} 
\label{eq:CDFSampleMean}\\
&=&
\frac{1}{{\mathfrak N}^N}
\sum\limits_{n_1,\dots,n_N=0}^{M} \left(\prod\limits_{p=1}^N \frac{a_{n_p}}{n_p!}\right)
\sum\limits_{l=(z+\mu) N}^{N M}
\sum\limits_{\sum\limits_{p=1}^N j_p = l}
\left(\prod\limits_{p=1}^N j_p^{n_p} 1_{j_p \le M}\right)
\label{eq:CDFSampleMean1}\\
&=&
\frac{1}{{\mathfrak N}^N}
\sum\limits_{n_1,\dots,n_N=0}^{M}
\frac{\prod\limits_{p=1}^N a_{n_p}}{(\left|\vec{n}\right|+N)!}
\left(
(NM)^{\left|\vec{n}\right|+N}
-
(z+\mu)^{\left|\vec{n}\right|+N}
\right)
\label{eq:CDFSampleMean2}\\
&=&
\frac{1}{{\mathfrak N}^N}
\frac{1}{(N-1)!}
\sum\limits_{p=0}^{{M} N}
\frac{(N M)^{p+N} - ((z+\mu)N)^{p+N}}{(p+N) p!} a^{\otimes N}_p
\label{eq:CDFSampleMean3} \\
&=&
\frac{1}{{\mathfrak N}^N}
\sum\limits_{l=(z+\mu)N}^{M N}
\frac{l^{N-1}}{(N-1)!}
\sum\limits_{p=0}^{{M} N}
\frac{l^{p}}{p!} a^{\otimes N}_p
\label{eq:CDFSampleMean4}
\end{eqnarray}
In (\ref{eq:CDFSampleMean}) we used the definition of the CDF and 
the Law of Conditional Probabilities, and we conditioned on the random variables in the sample, while we used equation (\ref{eq:ProbSpace}) in (\ref{eq:CDFSampleMean1}).
In (\ref{eq:CDFSampleMean2}) we neglected the indicator functions in the last term in parentheses, which we can do for large $M$, and we  summed over the $j$ values using  the identity (\ref{eq:SumovSimplexI}) and we summed over $l$ using  elementary number theoretic identities. 
We also defined 
$\left|\vec{n}\right| := \sum\limits_{p=1}^N n_p$
as the $L^1$ norm of the $\vec{n}$ vector.
In (\ref{eq:CDFSampleMean3}) we introduced the normalized $N$\textsuperscript{th} auto-convolution $a^{\otimes N}_p$ of the coefficients 
$\left(a_n\right)_{n=0}^{M}$ via:
\begin{equation}
a^{\otimes N}_p := \frac{1}{C^{N-1+p}_p} \sum\limits_{\sum\limits_{q=1}^N n_q = p} \prod\limits_{q=1}^N a_{n_q}
\label{eq:NthAutoConv}
\end{equation}

\noindent{\bf Note:} Expression (\ref{eq:CDFSampleMean4}) uses an approximate identity (\ref{eq:SumovSimplexI}) and is therefore only a starting point for analyzing the scaling of the percentiles of the difference between the sample mean and the population mean.
In order to provide further insight into this problem one needs to compute the $N$\textsuperscript{th} autoconvolution
from (\ref{eq:NthAutoConv}), (\ref{eq:InverseRel}) and from the expression (\ref{eq:InverseVandermonde})
 for the inverse Vandermonde matrix. This will make it possible to uniquely classify distributions
according to the asymptotic behaviour of their percentiles. We will accomplish that goal in future work.

\section{Conclusions}

All finite discrete probability distributions may be represented by polynomials in the possible outcomes of their random variables, by inverting the relevant VanderMonde matrix. 
We used this fact to examine the convergence of certain sample properties of finite discrete distributions. 
In particular, we considered families of probability distributions that had different numbers of outcomes ($M$) but were otherwise similar. 
Such families of distributions can be represented as sequences whose terms are indexed by $M$. Their variances can also be described as sequences in $M$.

In this paper we considered the special case of integer outcomes $\lbrace 0, 1, \ldots M \rbrace$. 
We calculated an  expression for the variances of such families of distributions using large $M$ approximations, and then examined their behaviour as $M \to \infty$. 
We found conditions for \textit{divergence} of the variance with increasing $M$ and gave some examples. We discussed some necessary and sufficient conditions for distributions to exhibit specific types of divergence. We also gave examples of families of distributions that \textit{did not} satisfy these conditions for divergence and whose variances actually \textit{converged} as $M \to \infty$. Finally, we have provided an expression for the Cumulative Distribution Function (CDF) of the difference between the sample mean and the population mean. This expression can be used to analyze the asymptotic behaviour of the percentiles of the difference between the sample mean and the population mean as a function of the number of outcomes $M$.

In obtaining the results, we derived many mathematical identities and they are listed in the appendices. We also derived an expression for the inverse of the Vandermonde matrix used in this paper, which was previously unknown. We believe these formulae could be useful in many different fields.

We hope, in future work, to generalize our results to arbitrary probability distributions, to quantify better what distributions these results apply to, as well as examine the behaviour of the variances of these distributions, for an increasing number of outcomes, relative to the mean. 

In conclusion, this work should be relevant to situations where one must sample the distributions of finite discrete random variables. It highlights how adding new outcomes to such random variables (without any change in the overall form of the probability distribution) may still affect how the statistics converge, with increasing sample size, to the population properties. More importantly, we have proven that for certain finite discrete probability distributions  when $N$ is kept \textit{fixed}, the variance of the sample mean actually \textit{diverges} as $M \to \infty$. We termed this result as the \textit{``Law of Many Outcomes''}, or alternatively, the \textit{``Reverse of the Law of Large Numbers''}. 

\appendix

\newcommand{\appsection}[1]{\let\oldthesection\thesection\renewcommand{\thesection}{Appendix \oldthesection}\section{#1}\let\thesection\oldthesection}

\appsection{}\label{AppA} 

In this appendix, we list certain identities used in the main body of the paper. The identities (\ref{eq:Faulhaber})-(\ref{eq:SumSimplexV_A}) all relate to sums of powers of integers over certain sets. The identities (\ref{eq:BinomIdentity1})-(\ref{eq:BinomIdentity}) involve sums of binomial coefficients. Then the remainder of this appendix is a generalization of the standard proof of the Vandermonde determinant. This is used in Appendix B to calculate the inverse of the Vandermonde matrix used in this paper.

In the following we denote by $\left(B_k\right)_{k=0}^\infty$ the Bernoulli numbers.

\noindent{\bf Faulhaber's Formula:}

\begin{equation}
\sum\limits_{j=0}^M j^n = \sum\limits_{k=0}^n (M+1)^{n-k+1} \frac{n!}{(n-k+1)!} \frac{B_k}{k!}
\label{eq:Faulhaber}
\end{equation}

\noindent{\bf Sum Over Simplex I:}

\begin{equation}
\sum\limits_{\sum\limits_{p=1}^N j_p = l}
\prod\limits_{p=1}^N j_p^{n_p}
=
l^{\left|\vec{n}\right|+N-1}
\frac{\prod\limits_{p=1}^N n_p!}{(\left|\vec{n}\right|+N-1)!}
\label{eq:SumovSimplexI}
\end{equation}
Here $\left|\vec{n}\right| := \sum\limits_{p=1}^N n_p$.
The result is valid for big values of $\left|\vec{n}\right|$ only.

\noindent{\bf Sum Over Simplex II:}

\begin{equation}
\sum\limits_{0\le j_1<\dots<j_{s}\le n} 
\prod\limits_{p=1}^{s} j_p^{n_p}
=
\sum\limits_{m=1}^{N_s+s}
(n+1)^m
{\mathbb C}^{(s)}_m
\label{eq:SumovSimplexII}
\end{equation}
Here $N_j := \sum\limits_{p=1}^j n_p$ is the $L^1$-norm of the sequence $\left(n_p\right)_{p=1}^j$, the symbol $(a)_{(n)} := \prod\limits_{p=0}^{n-1} (a-p)$ is the Pochhammer symbol, and the coefficients ${\mathbb C}^{(s)}_m$ read:
\begin{eqnarray}
{\mathbb C}^{(s)}_{N_s+s}   &:=& \frac{1}{\prod\limits_{q=1}^s (N_q+q)} \\ 
{\mathbb C}^{(s)}_{N_s+s-1} &:=& \sum\limits_{j_1=1}^s \frac{(N_{j_1}+j_1)}{\prod\limits_{q=1}^s (N_q+q-1_{q>j_1})} \frac{B_1}{1!} \\
{\mathbb C}^{(s)}_{N_s+s-2} &:=& \sum\limits_{j_1=1}^s \frac{(N_{j_1}+j_1)_{(2)}}{\prod\limits_{q=1}^s (N_q+q-2\cdot 1_{q>j_1})} \frac{B_2}{2!} \nonumber \\
&& {}+ \sum\limits_{1\le j_1 < j_2 \le s} \frac{(N_{j_1}+j_1)_{(1)}(N_{j_2}+j_2-1)_{(1)}}{\prod\limits_{q=1}^s (N_q+q-1_{q>j_1}-1_{q>j_2})} (\frac{B_1}{1!})^2 \\
{\mathbb C}^{(s)}_m &:=& 
\sum\limits_{1\le j_1\le \dots \le j_{N_s+s-m} \le s}
\frac{\prod\limits_{p=1}^{N_s+s-m}(N_{j_p}+j_p+1-p)}{\prod\limits_{p=1}^s(N_p+p-\sum\limits_{l=1}^x 1_{p>j_l})}
\cdot \prod\limits_{p=1}^{l} \frac{B_{d_p}}{d_p!}
\label{eq:SumovSimplexIII}
\end{eqnarray}
for $m=1,\dots,N_s+s$.
Here the numbers $d_p$ are multiplicities of elements of the sequence $\left(j_p\right)_{p=1}^s$, 
i.e.\ such numbers that the sequence $\left(j_p\right)_{p=1}^s$ falls into $l$ groups composed of equal elements, such that the first group has length $d_1$,
the second group length $d_2$, up to the $l$\textsuperscript{th} group who has length $d_l$.

In particular when $n_j = 1$  and thus $N_j = j$ for $j=1,\dots,s$ then we have:
\begin{equation}
{\mathbb C}^s_{m} := \frac{1}{m!!}
\sum\limits_{1\le j_1 \le \dots \le j_{2s-m} \le s}
\prod\limits_{p=1}^{2s-m} \frac{(2 j_p - p)!!}{(2 j_p - (p+1))!!}
\prod\limits_{p=1}^l \frac{B_{d_p}}{d_p!}
\label{eq:SumovSimplexSpecial}
\end{equation}

We were able to resum the series for particular values of $m$, using identity (\ref{eq:HypergeometricIden}), and we give the result below:
\begin{eqnarray}
{\mathbb C}^{(s)}_{2s}   &=& \frac{1}{(2s)!!} \nonumber \\
{\mathbb C}^{(s)}_{2s-1} &=& \frac{1}{(2s-1)!!} \left( \frac{1}{3}\frac{(2s+1)!!}{(2s-2)!!}  B_1 \right)\nonumber\\
{\mathbb C}^{(s)}_{2s-2} &=& \frac{1}{(2s-2)!!} \left( \frac{1}{3^2} (s-1) s (1+4s)B_1^2 + s^2 \frac{B_2}{2!} \right) \\
{\mathbb C}^{(s)}_{2s-3} &=& \frac{1}{(2s-3)!!} \frac{(2s-1)!!}{(2s-4)!!}\left( 
 \frac{1}{3} \left(\frac{2}{45} + \frac{7}{135} s - \frac{1}{3} s^2 + \frac{4}{27} s^3\right) B_1^3  \right. \nonumber \\
&& \left. + \frac{1}{3} \left(-\frac{1}{5}-\frac{2}{5} s + s^2\right) \frac{B_2}{2!} B_1 + \frac{1}{5} (2s+1)\frac{B_3}{3!} \right) \\
\vdots \nonumber
\label{eq:SumovSimplexSpecialI}
\end{eqnarray}

The identity (\ref{eq:SumovSimplexII}) is proven by performing the sums over the consecutive indices $j_{s},\dots,j_1$ by means of 
the Faulhaber's formula.

\noindent{\bf Sum Over Simplex III:}

\begin{eqnarray}
\sum\limits_{a\le j_1<\dots<j_{s}\le b} 
\prod\limits_{p=1}^{s} j_p^{n_p}
=
\mathop{
\sum\limits_{m_1=0}^{N_s+s}
\sum\limits_{m_2=0}^{N_s}
}_{m_1+m_2\le N_s+s}
{\mathbb C}^{(s)}_{m_1} C^{N_s}_{m_2} (b-a+1)^{m_1} a^{m_2}
\label{eq:SumovSimplexIV}
\end{eqnarray}
where the coefficients ${\mathbb C}^{(s)}_{m}$ are defined in (\ref{eq:SumovSimplexIII}).
The identity (\ref{eq:SumovSimplexIV}) follows from substituting $j^{'}_q = j_q - a$ for $q=1,\dots,s$, from expanding the
resulting power terms in binomial expansions, from applying  (\ref{eq:SumovSimplexII}) to sum over the sequences 
$\left(j^{'}_q\right)_{p=1}^s$ and from the identity (\ref{eq:BinomIdentity}).

\noindent{\bf Sum Over Simplex IV:}

\begin{equation}
{\mathfrak S}^{\tilde{q}}_{q,j} :=
\sum\limits_{q< j_{\tilde{q}-1} < \dots < j_0 < j} \prod\limits_{l=0}^{\tilde{q}} C^{j_{l-1}-1}_{j_l-1}
=
C^{j-1}_{q-1} \sum\limits_{l=1}^{\tilde{q}+1} l^{j-q} (-1)^{\tilde{q}+1-l} C^{\tilde{q}+1}_{\tilde{q}+1-l}
\label{eq:SumSimplexIV}
\end{equation}
Here $j_{-1} = j$ and $j_{\tilde{q}}= q$ and $j-q \ge \tilde{q}+1$.
The series (\ref{eq:SumSimplexIV}) can be resummed by means of the binomial expansion formula.

As a simple corollary from (\ref{eq:SumSimplexIV}) we have:

\noindent{\bf Sum Over Simplex V:}

\begin{eqnarray}
\lefteqn{
\sum\limits_{q+p_1 < j_{\tilde{q}-1} < \dots < j_0 < j - p_2} \prod\limits_{l=0}^{\tilde{q}} C^{j_{l-1}-1}_{j_l-1}
=
{\mathfrak S}^{\tilde{q}}_{q,j}
-
\nonumber }\\
&&\!\!\!\!\!\!\!\!\! \mathop{\mathop{\sum\limits_{l_1=1,\dots,p_1}}_{l_2=1,\dots,j-q-(p_1+p_2)}}_{l_3=1,\dots,p_2+1}
\!\!\!\!\!\!\!\!\!\!\!\!\!\!\!\!\!\!
C^{q+p_1-1}_{q-1} C^{j-p_2-1}_{q+p_1-2} C^{j-1}_{j-p_2-2}
l_1^{p_1} l_2^{j-q-(p_1+p_2)} l_3^{p_2+1}
(-1)^{\tilde{q}+3-(l_1+l_2+l_3)}
\beta_{l_1,l_2,l_3}^{p_1,p_2}
\label{eq:SumSimplexV}
\end{eqnarray}
where 
\begin{equation}
\beta_{l_1,l_2,l_3}^{p_1,p_2}(\tilde{q}) := 
\sum\limits_{0 < q_2 < q_1 < \tilde{q}}
C^{\tilde{q}-1-q_1+1}_{l_1}
C^{q_1-q_2+1}_{l_2}
C^{q_2+1}_{l_3}
1_{q_2 \le p_2} 1_{q_1-q_2 \le j-q-(p_1+p_2)} 1_{\tilde{q}-q_1 \le p_1-1}
\label{eq:SumSimplexV_A}
\end{equation}
Here the second bit on the right hand side in (\ref{eq:SumSimplexV}) stands for the sum over ordered sequences
bounded from below and above by $q$ and $j$ respectively and not bounded from below and above by $q+p_1$ and $j-p_2$.
Here $j-q \ge p_1+p_2 + \tilde{q}+1$.

\noindent{\bf A Binomial Identity:} 

Let $M> s \in {\mathbb N}$ and $p\in {\mathbb Z}$ and $p \ge -1$. Then we have:
\begin{eqnarray}
\lefteqn{ 
\sum\limits_{j=1}^{M-s} C^{x+j-1}_x (j+s)^{p+1}
}  \nonumber \\
&&
= C^{M-s+x}_{x+1} \left((x+1)(M-s-1)!\sum\limits_{q=1}^{M-s} \frac{(x+q-1)!}{(q-1)! (x+M-s)!(s+q)^{p+1}} \right) 
\label{eq:BinomIdentity1} \\
&&=C^{M-s+x}_{x+1} f_p(M-s,x) \label{eq:BinomIdentity2}
\end{eqnarray}
Note that for $p=-1$ the term in parentheses on the right hand side equals unity.
The quantities $f_p(n,x)= (x+1)\sum\limits_{q=0}^{p+1} A_q^{p+1} n^q$ are polynomials of order $p+1$  
in the variable $n$ with coefficients $A_q^{p+1} = A_q^{p+1}(s)$ that depend on $s$.
Following recursion relations hold
\begin{equation}
 f_p(n,x) = \frac{(n-1)}{(x+n)} f_p(n-1,x) + \frac{(x+1)}{(x+n)} (s+n)^{p+1}
= \left( 1, \frac{s(x+2)+1}{(x+1)(x+2)} + \frac{x+1}{x+2} n, \dots \right)_{p=-1}^\infty
\end{equation}
for the polynomials and
\begin{equation}
A_{q-1}^{p+1} 1_{q\ge 2} + x A_q^{p+1} 1_{q\le p+1}= \sum\limits_{l=q-1}^{p+1} A_l^{p+1} C^{l+1}_q (-1)^{l+1-q} + C^{p+1}_q s^{p+1-q} 1_{q\le p+1}
\label{eq:PolynomCoeffs}
\end{equation}
and for the coefficients of the polynomials. Here $q=p+1,p,\dots,0$.

Below we enclose a Mathematica script that test the identities (\ref{eq:BinomIdentity1}) and (\ref{eq:BinomIdentity2}):

\textit{(* Testing the binomial identity (59)  *)}

\begin{verbatim}
s=1;p=6;x=7;

Table[Sum[ Binomial[ x+j-1,x](j+s)^{}(p+1),{j,1,M-s}],{M,2,20}]
-
Table[Binomial[M-s+x,x+1] Sum[(M-s-1)!/(q-1)! (x+q-1)!/(x+M-s)! (s+q)^{}(p+1),{q,1,M-s}](x+1),{M,2,20}]
\end{verbatim}

\textit{(* Testing the binomial identity (60)  *)}

\begin{verbatim}
x=.;n=.;s=.;p=3;

ll=CoefficientList[ (x+n)Sum[A[q] n^{}q,{q,0,p+1}]- (n-1)Sum[A[q] (n-1)^{}q,{q,0,p+1}] - (s+n)^{}(p+1),{n}];

ss =Simplify[Solve[Table[ ll[[j]] == 0,{j,1,Length[ll]}],Table[A[q],{q,0,p+1}]]];

MyCoeffs =Table[ A[q] ,{q,0,p+1}] /. ss[[1]]

MyF[n_,x_,s_] := Sum[ MyCoeffs[[i]] n^{}(i-1),{i,1,Length[MyCoeffs]}];

s=4;x=7;

Table[Sum[Binomial[ x+j-1,x](j+s)^{}(p+1),{j,1,M-s}],{M,5,30}]-

Table[Binomial[M-s+x,x+1]  (x+1) MyF[M-s,x,s],{M,5,30}]
\end{verbatim}

\noindent{\bf A Hypergeometric Identity:}

Let $a\in {\mathbb N}$ and $a_q \in {\mathbb R}$ for $q=0,\dots,N$. Then we have:
\begin{equation}
\sum\limits_{j=1}^s \frac{(2j + 2a-1)!!}{(2j-2)!!} (\sum\limits_{q=0}^N a_q j^q) = 
\frac{1}{2a+3} \frac{(2s+2a+1)!!}{(2s-2)!!} (\sum\limits_{q=0}^N A_q s^q)
\label{eq:HypergeometricIden}
\end{equation}
where the coefficients $\left(A_p\right)_{p=0}^N$ satisfy the following system of equations:
\begin{equation}
\left((2a+1)+2(p+1)\right) A_p - 2 \sum\limits_{q=p+1}^N C^{q+1}_p (-1)^{q+1-p} A_q = (2a+3) a_p
\end{equation}
for $p=N,N-1,\dots,0$. The system of equations always has a unique solution.
The identity (\ref{eq:HypergeometricIden}) is proven by elementary methods, i.e.\ by reducing the sum to a common denominator and factoring out the numerator. 
It would be interesting to know if an analoguous identity exists in the case when the ratio of double factorials is replaced by a product of one or several ratios of that kind or of similar ones.

\noindent{\bf Corollary}: 
Define:
\begin{equation}
{\mathfrak S}^{(d)}(s) := 
\sum\limits_{1\le j_1 < \dots < j_d \le s}
\prod\limits_{p=1}^d \frac{(2j_p-p)!!}{(2j_p-(p+1))!!}
\label{eq:simplBulk}
\end{equation}
and
\begin{equation}
{\mathfrak S}^{(d,i)}(s) := 
\sum\limits_{1\le j_1 < \dots < j_d \le s}
\prod\limits_{p=1}^d \frac{(2j_p-p)!!}{(2j_p-(p+1))!!}\delta_{j_i,j_{i+1}}
\label{eq:simplBorder}
\end{equation}
for $i=1,\dots,d-1$. We term the quantities in (\ref{eq:simplBulk}) and in (\ref{eq:simplBorder}) the sum over the bulk
and over the border of a $d$ dimensional simplex, respectively.

Then the sums over the bulk of the simplex read:
\begin{eqnarray}
{\mathfrak S}^{(2)}(s)
&=&
\frac{1}{9}s(s-1)(4s+1) \nonumber \\
{\mathfrak S}^{(4)}(s)
&=&
\frac{1}{2430}(s-1)(s-2)(s-3)(80 s^3-120 s^2+7 s + 12) \nonumber \\
{\mathfrak S}^{(6)}(s)
&=&
\frac{1}{2296350}(s-2)(s-3)(s-4)(s-5) (2240 s^5-14000s^4+27580 s^3-17815 s^2+159 s + 1566) \nonumber \\
{\mathfrak S}^{(3)}(s)
&=&
\frac{1}{405}\frac{(2s-1)!!}{(2s-4)!!} (s-2) \left(20 s^2-5 s -3\right) \nonumber \\
{\mathfrak S}^{(5)}(s)
&=&
\frac{1}{51030}
\frac{(2s-3)!!}{(2s-6)!!}
(s-3)(s-4) (112 s^4 - 392 s^3 + 329 s^2 - 7s -30) \nonumber \\
{\mathfrak S}^{(7)}(s)
&=&
\frac{1}{6889050}
\frac{(2s-5)!!}{(2s-8)!!} \times
\nonumber \\
&& (s-4)(s-5)(s-6)(-756 - 15 s + 8894 s^2 - 16035 s^3 + 10820 s^4 - 3120 s^5 + 320 s^6)  
\label{eq:SumBulkExpl}
\end{eqnarray}

The sums over the border of the simplex read:
\begin{eqnarray}
{\mathfrak S}^{(3,1)}(s)
&=&
\frac{(2s-1)!!}{(2s-4)!!} \frac{1}{3} \left(-\frac{1}{35}-\frac{12}{35}s + \frac{3}{7} s^2\right) \nonumber \\
{\mathfrak S}^{(3,2)}(s)
&=&
\frac{(2s+1)!!}{(2s-4)!!} \frac{1}{5} \left(-\frac{6}{21}+\frac{10}{21}s\right) \nonumber \\
{\mathfrak S}^{(4,1)}(s)
&=&
\frac{1}{630}(s-1)(s-2)(3s-2)(24 s^2-41 s-3) \nonumber \\
{\mathfrak S}^{(4,2)}(s)
&=&
\frac{2}{105}(s-1)(s-2)(2s-3)(2s-1)(2s+1) \nonumber \\
{\mathfrak S}^{(4,3)}(s)
&=&
\frac{1}{90}(s-1)(s-2)s(16 s^2-17s-3) 
\label{eq:SumBorderExpl}
\end{eqnarray}
The identities (\ref{eq:SumBulkExpl}) and (\ref{eq:SumBorderExpl}) follow from an iterative application of the 
identity (\ref{eq:HypergeometricIden}). The sums are performed 
starting from the sum over the index with the lowest subscript and ending at the sum over the index with the biggest subscript.
It will be interesting to find a generic formula related to a simplex of arbitrary dimension.
In future work we will use the above results to conjecture the generic formula and prove that formula by induction.

\noindent{\bf The Generalized Chu-Vandermonde Identity:}

\begin{equation}
\sum\limits_{0\le q_1 \le \dots \le q_{s-1} \le b} \prod\limits_{p=1}^s C^{a_p}_{q_p-q_{p-1}} =
C^{\sum\limits_{p=1}^s a_p}_{b}
\label{eq:BinomIdentity}
\end{equation}
for $q_0=0$ and $q_s = b$.
The above identity is proven by an iterative application of the Chu-Vandermonde identity
\cite{ChuVandermonde}.
Since the Chu-Vandermonde identity is closely related to the Gauss's hypergeometric theorem \cite{GaussHypergeomTh}, to Dougall's Formula \cite{Dougall}, to Thomae's theorem \cite{Thomae} and to various other identities that involve generalized hypergeometric functions \cite{Ramanujan}, it would be interesting to derive the many-dimensional analogues for those identities using our methods.

\noindent{\bf Generalized Vandermonde Determinant:}

Let $M\in {\mathbb N}$ and $1 \le a \le M$ and 
$1\le J_1< \dots < J_a \le M$  with $J_0=0$ and $J_{a+1}=M$, and  $\delta_\theta \in {\mathbb N}$ for $\theta=1,\dots,a$. 
In addition define $\Delta_j := \sum\limits_{\theta=1}^j \delta_\theta$, and
$\vec{x}_a^b := \left(x_j\right)_{j=a}^b$ and
\begin{equation}
P^{(l)}\left(\vec{x}_a^b\right) := \sum\limits_{a\le j_1<\dots<j_l\le b} \prod\limits_{q=1}^l x_{j_q} 
\end{equation}

Then we have:

\begin{eqnarray}
\lefteqn{\det\left( x_p^{j-1 + \sum\limits_{\theta=1}^a \delta_\theta 1_{j\ge J_\theta+1}} \right)_{p,j=1,1}^{M,M}=}
\label{eq:GenVandeMondeFull0} \\
&&
=\left(\prod\limits_{M\ge j>1} (x_j-x_1)\right) \cdot
\sum\limits_{\xi_1,\dots,\xi_a=0}^{\delta_1,\dots,\delta_a}
x_1^{\sum\limits_{\theta=1}^a (\delta_\theta - \xi_\theta)}
\det\left( x_p^{j-2 + 
\sum\limits_{\theta=1}^a (\delta_\theta 1_{j\ge J_\theta+1}
+ (\xi_\theta - \delta_\theta) 1_{j=J_\theta+1})} \right)_{p,j=2,2}^{M,M}
\label{eq:GenVandeMondeFull}\\
&&=\left(\prod\limits_{M\ge j>1} (x_j-x_1)\right) \cdot
\sum\limits_{\vec{\xi} \in \otimes_{\theta=1}^a [0,\delta_\theta]}
x_1^{\left|\vec{\delta}-\vec{\xi}\right|} 
\det\left( x_{p+1}^{j-1 + 
\sum\limits_{\theta=1}^a
( \xi_\theta 1_{j\ge J_\theta} + (\delta_\theta - \xi_\theta)1_{j\ge J_\theta+1} 1_{J_{\theta+1}-J_\theta \ge 2})
}
 \right)_{p,j=1,1}^{M-1,M-1}
\label{eq:GenVandeMondeFull1}\\
&&=
\left(\prod\limits_{M\ge j>i \ge 1} (x_j-x_i)\right) \cdot
\sum\limits_{\vec{\xi}^{1} \in \otimes_{\theta^1=1}^{a^1} [0,\delta^{1}_{\theta^1}]}
\cdot
\dots
\cdot
\sum\limits_{\vec{\xi}^{M} \in \otimes_{\theta^M=1}^{a^M} [0,\delta^{a^M}_{\theta^M}]}
\prod\limits_{p=1}^M x_p^{\left|\vec{\delta}^{p}-\vec{\xi}^{p}\right|}
\label{eq:GenVandeMondeFull2}
\end{eqnarray}
The sum on the right-hand-side in (\ref{eq:GenVandeMondeFull2}) contains 
$\prod\limits_{j=1}^M \prod\limits_{q_j=1}^{a^j} (\delta^j_{q_j}+1)$ terms.
Here the parameters $\left(a^{p+1},\vec{J}^{p+1},\vec{\delta}^{p+1}\right)_{p=0}^{M-1}$ 
with $\left|\vec{J}^{p+1}\right| = \left|\vec{\delta}^{p+1}\right| = a^{p+1}$
for $p=0,\dots,M-1$
constitute a branching process
and satisfy following recursion relations:
\begin{eqnarray}
a^{p+1} &:=& a^p  + \sum\limits_{\theta=1}^{a^p} 1_{J^p_{\theta+1} - J^p_{\theta} \ge 2} \nonumber\\
\vec{J}^{p+1} &:=& \left(\underbrace{J_1^p-1,\dots,J_{A_1}^p-1}_{A_q},
\underbrace{J_{A_1+1}^p-1,J_{A_1+1}^p,\dots,J_{A_1+A_2}^p-1,J_{A_1+A_2}^p}_{2 A_{q+1}}\right)_{q=1}^\omega\nonumber\\
\vec{\delta}^{p+1} &:=& \left(\underbrace{\xi_1^p,\dots,\xi_{A_1}^p}_{A_q},
\underbrace{\xi_{A_1+1}^p,\delta_{A_1+1}^p-\xi_{A_1+1}^p,\dots,\xi_{A_1+A_2}^p,\delta_{A_1+A_2}^p-\xi_{A_1+A_2}^p}_{2 A_{q+1}}\right)_{q=1}^\omega
\label{eq:RecursRel}
\end{eqnarray}
for $p=1,\dots,M-1$.
The constraint, $\sum\limits_{q=1}^\omega (A_q + 2A_{q+1}) = a^{p+1}$, holds where $A^q$ counts the number of consecutive adjacent elements in the sequence $\vec{J}^p$ and $A_{q+1}$ counts the number of consecutive non-adjacent elements in the sequence $\vec{J}^p$.
Here $q=1,\dots,\omega$.

\noindent{\bf Proof:} 
We denote by $Q^{(n)}(x_1,x_2) := \mathop{\sum\limits_{l_1+l_2=n}}_{l_1,l_2\ge 0} x_1^{l_1} x_2^{l_2}$ a totally symmetric polynomial of order $n$
in the variables $x_1$ and $x_2$ and by $n_j:=j-2+\sum\limits_{\theta=1}^a \delta_\theta 1_{j\ge J_\theta+1}$ for $j=2,\dots,M$
the power index of the matrix element.
Now, we multiply the first row of the determinant by minus unity and add to all following rows, i.e.\ to the second, the third, up to
the $M$\textsuperscript{th} row, we factor out a term $(x_p-x_1)$ from the $p$\textsuperscript{th} row, and, from the multi-linearity of the determinant, we get:
\begin{equation}
\det = \left(\prod\limits_{p=2}^M (x_p-x_1)\right) \det\left(Q^{(n_j)}(x_1,x_p)\right)_{p,j=2,2}^{M,M}
\end{equation}
In the next step we modify the $j$\textsuperscript{th} column, for $j=2,\dots,M$.  
We multiply the first column by $(-x_1^{n_j})$, the second column by $(-x_1^{n_j-1})$,
and so on and so forth, up to the $J_1$\textsuperscript{th} column by $(-x_1^{n_j-J_1+1})$ and add them all to the $j$\textsuperscript{th} column.
It is not hard to see that we get:
\begin{eqnarray}
\lefteqn{\det\left(\sum\limits_{l=J_1-1}^{n_j} x_p^{l} x_1^{n_j-l})\right)_{p,j=2,2}^{M,M}
=}
\nonumber \\
&&\det\left(1,x_p,x_p^2,\dots,x_p^{J_1-2}, \sum\limits_{l_1=J_1-1}^{n_j} x_p^{l_1} x_1^{n_j-l_1}, \dots, 
\sum\limits_{l_{M-J_1}=J_1-1}^{n_j} x_p^{l_{M-J_1}} x_1^{n_j-l_{M-J_1}}\right)_{p=2}^M
\end{eqnarray}
Here in the first $J_1-1$ columns (whose labels are $2,\dots,J_1$) we have single powers terms,
whereas in the last $M-J_1$ columns (whose labels are $J_1+1,\dots,M$) we have sums of $n_j-J_1+2$ power terms.
Thus at positions $j=J_\theta+1$ for $\theta=1,\dots,a$ the number of terms in the sums increases by 
$\delta_1,\delta_2+1,\dots,\delta_a+1$
whereas at all remaining positions, with $j\ge J_1+1$, the number of terms in the sums increases by one. 
Now we expand the determinant and we deduce, from the multi-linearity and antisymmetric property of the determinant that the sums at positions $j=J_\theta+1$
can run over the whole range of index values, 
whereas  the sums at the remaining positions pick up only maximal values of the summation index.
Otherwise in the determinants that result from the expansion  there are at least two columns that are proportional
to each other, hence the determinants are zero. In other words we have:
\begin{eqnarray}
l_1 &\in & \left\{J_1-1,\dots,J_1-1+\delta_1\right\},  l_2 = J_1+\delta_1,\dots, l_{J_2-J_1} = J_2-2+\delta_1 \\
l_{J_2-J_1+1} &\in & \left\{J_2-1+\delta_1,\dots,J_2-1+\Delta_2\right\}, l_{J_2-J_1+2} = J_2+\Delta_2,\dots, 
\nonumber \\ 
&& l_{J_3-J_2} = J_3-2+\Delta_2\\
&&\vdots \nonumber \\
l_{J_a-J_{a-1}+1} &\in& \left\{J_a-1+\Delta_{a-1},\dots,J_a-1+\Delta_a\right\}, l_{J_a-J_{a-1}+2} = J_a+\Delta_a,\dots,  \nonumber \\ 
&& l_{M-J_1} = M-2+\Delta_a
\end{eqnarray}
From this follows the result (\ref{eq:GenVandeMondeFull}).
The final result (\ref{eq:GenVandeMondeFull2}) follows from iterating (\ref{eq:GenVandeMondeFull1}).
{\bf q.e.d}.

{\bf Example 1:} Take $a=1$ and $\delta_1 = 1$. Then we have:
\begin{eqnarray}
\lefteqn{
\det\left(x_p^{j-1+1_{j\ge J+1}}\right)_{p,j=1,1}^{M,M} = }
\nonumber\\
&&
\left(\prod\limits_{M\ge j>1} (x_j-x_1)\right) \cdot
\left[
x_1 \det\left(x_p^{j-1+1_{j\ge J+1}}\right)_{p,j=2,1}^{M,M-1} +
\det\left(x_p^{j-1+1_{j\ge J}}\right)_{p,j=2,1}^{M,M-1}
\right]=
\label{eq:GenVandeMonde1} \\
&&
\left(\prod\limits_{M\ge j>i\ge 1} (x_j-x_i)\right) \cdot
\sum\limits_{1\le i_1 < \dots < i_{M-J} \le M} \prod\limits_{q=1}^{M-J} x_{i_q}
\label{eq:GenVandeMonde2} 
\end{eqnarray}
The last step (\ref{eq:GenVandeMonde2}) is proven by mathematical induction in $M$, for example.

\noindent{\bf Example 2:} Take $a=1$ and $\delta_1>0$. Then we have:
\begin{eqnarray}
\lefteqn{
\det\left(x_p^{j-1+\delta_1 1_{j\ge J+1}}\right)_{p,j=1,1}^{M,M} = \left(\prod\limits_{M\ge j>i\ge 1} (x_j-x_i)\right) \cdot}
\nonumber\\
&&
\sum\limits_{0\le \xi_1^M \le \dots \le  \xi_1^1 \le \delta_1}
\sum\limits_{0\le \xi_2^M \le \dots \le  \xi_2^2 \le \delta_1-\xi_1^1}
\cdots
\sum\limits_{0\le \xi_{M}^M \le \delta_1-\sum\limits_{q=1}^{M-1} \xi_q^q}
\prod\limits_{q=1}^M x_q^{ 
\sum\limits_{q=1}^{M-1} \xi^{M-1}_q
+\delta_1 - (\sum\limits_{q=1}^{M-1} \xi^q_q)
-(\sum\limits_{q=1}^{M} \xi^{M}_q)
}
\label{eq:GenVandeMonde3} 
\end{eqnarray}
The result (\ref{eq:GenVandeMonde3}) follows from (\ref{eq:GenVandeMondeFull2})
along with the recursion relations (\ref{eq:RecursRel}).

\appsection{}\label{AppB} 

Here we prove formula (\ref{eq:InverseVandermonde}) for the inverse of the Vandermonde matrix
\mbox{$\left({\mathfrak A}_{j,n}\right)_{j=0,n=1}^{M-1,M} = \left[\left(j^{n-1}\right)_{j,n=1}^M\right]^{-1}$}.
\begin{eqnarray}
\lefteqn{
{\mathfrak A}_{j,n} = \frac{(-1)^{j+n+1}}{\prod\limits_{1\le p < q \le M} (q-p)}
\det\left((p+1_{p\ge j})^{q-1+1_{q\ge n+1}}\right)_{p,q=1}^{M-1}
}\label{eq:MatrixInverse}\\
&&= 
\left( \frac{(-1)^{j+n+1}\prod\limits_{1\le p < q \le M-1} (q-p + 1_{q\ge n} - 1_{p\ge n})}{\prod\limits_{1\le p < q \le M} (q-p)} \right) \times \nonumber \\
&& \sum\limits_{1\le i_1 < \dots < i_{M-1-j} \le M-1}
\prod\limits_{q=1}^{M-1-j} (i_q+1_{i_q\ge n})
\label{eq:MatrixInverse1} \\
&&=
\frac{(-1)^{j+n+1}}{(n-1)!(M-n)!}
\cdot
S^{M-1-j}(M-1)
\label{eq:MatrixInverse1a}
\end{eqnarray}

In (\ref{eq:MatrixInverse}) we expressed the inverse matrix as the transposed matrix of algebraic complements and in 
(\ref{eq:MatrixInverse1}) we used
identity (\ref{eq:GenVandeMonde2}) to compute the algebraic complements.
In (\ref{eq:MatrixInverse1a}) we evaluated the sum of the product in a recursive way as follows:
\begin{eqnarray}
S^{M-1-j}(M-1) &:=& 
\sum\limits_{1\le i_1 < \dots < i_{M-1-j} \le M-1}
\prod\limits_{q=1}^{M-1-j} (i_q+1_{i_q\ge n})
=
\sum\limits_{l=1}^{M-1} S^{M-2-j}(l-1) (l + 1_{l\ge n})\nonumber \\
&=&
\left(\sum\limits_{p=0}^{M-1-j} P_{j+1+p}(M) (-n)^p\right)
\end{eqnarray}
with $S^{0}(M-1) = 1$. The last equality has been obtained by computing the sums $S^{M-1-j}(M-1)$ 
recursively for $j=M-1,M-2,M-3,\dots$. Here we noticed that the result is an order-$j$ polynomial in the column number $n$,
with the coefficient at the $p$\textsuperscript{th} power of $n$ being a polynomials in $M$ of order $2j - 2p$. 
The later polynomials satisfy recursion relations (\ref{eq:PolynomialsRecurs}), which we verified using Mathematica
using the piece of code provided in Appendix C.
Changing the variables $\tilde{j} = j+1$ and $j = \tilde{j}$ yields expression (\ref{eq:InverseVandermonde}).

Now, we prove formula  (\ref{eq:InverseVandermonde1}).
In (\ref{eq:MatrixInverse1}) we write \mbox{$i_q = q + \sum\limits_{p=1}^j 1_{q\ge q_p-(p-1)}$} for $q=1,\dots,M-1-j$ and for some $1\le q_1 < q_2 < \dots < q_j \le M-1$.
We split the product under the sum into a product of $(j+1)$ products such that, in the $p$\textsuperscript{th} product, the index $q$ runs within the range $q=q_p-(p-1),\dots,q_{p+1}-(p+1)$ for $p=0,\dots,j$ with $q_{j+1}=M-1$.
Then we change variables $q_p^{'} = q+p$ for $q=q_p-(p-1),\dots,q_{p+1}-(p+1)$.
Then we assume that the column index $n$ satisfies the following inequality
$q_p+1 \le n \le q_{p+1}$. 
This implies that in the first $(p-1)$ products (from left to right), all the indicator functions equal zero, in the $p$\textsuperscript{th} product some indicator functions are zero and the others are one, and finally, in the last $(p+1)$ products all the indicator functions equal unity. 
This allows us to factor out the term $M!/n$ from the product and absorb it into the prefactor in (\ref{eq:MatrixInverse1}) thus producing the binomial coefficient. 
The remaining sum over the values of the $q$ indices is left unevaluated. 
It is possible to obtain the large $M$ limit of that sum easily.

\appsection{}\label{AppC} 

Here we give closed form expressions for certain polynomials that
are used to invert the Vandermonde matrix. The expressions have been obtained by solving
the recursion relations in (\ref{eq:PolynomialsRecurs}).
Define 
\mbox{$\tilde{P}_{M-j}(M) := P_{M-j}(M)/\left(\prod\limits_{p=-1}^{j-1} (M-p)\right)$}. Then we have:
\begin{eqnarray}
\tilde{P}_{M-1}(M) &=& \frac{1}{2} \label{eq:Pol1}\\   
\tilde{P}_{M-2}(M) &=& \frac{1}{24}  (3 M+2)\label{eq:Pol2}\\   
\tilde{P}_{M-3}(M) &=& \frac{1}{48} M (M+1) \label{eq:Pol3}\\
\tilde{P}_{M-4}(M) &=& \frac{ \left(15 M^3+15
    M^2-10 M-8\right)}{5760}\label{eq:Pol4}\\
\tilde{P}_{M-5}(M) &=& \frac{M (M+1) \left(3
    M^2-M-6\right)}{11520}\label{eq:Pol5}\\
\tilde{P}_{M-6}(M) &=& \frac{\left(63
    M^5-315 M^3-224 M^2+140 M+96\right)}{2903040}\label{eq:Pol6}\\
\tilde{P}_{M-7}(M) &=& \frac{M (M+1)
    \left(9 M^4-18 M^3-57 M^2+34 M+80\right)}{5806080}\label{eq:Pol7}\\
\tilde{P}_{M-8}(M) &=& \frac{\left(135 M^7-315 M^6-1575 M^5+735 M^4+5320 M^3\right)}{1393459200}\nonumber \\
&& \mbox{} + \frac{\left(2820 M^2-1936 M-1152\right)}{1393459200}\label{eq:Pol8}\\
\tilde{P}_{M-9}(M) &=& \frac{M (M+1) \left(15 M^6-75 M^5-135 M^4+527 M^3+768 M^2 \right)}{2786918400} \nonumber \\
&& \mbox{}- \frac{M (M+1) \left(668 M +1008\right)}{2786918400}\label{eq:Pol9}\\
\tilde{P}_{M-10}(M) &=& \frac{
     \left(99 M^9-594 M^8-1386 M^7+6468 M^6+14091 M^5\right)}{367873228800} \nonumber \\
&& \mbox{}- \frac{\left(12826 M^4 + 44132 M^3 + 18392 M^2 - 14432 M - 7680\right)}{367873228800}\label{eq:Pol10}\\
\vdots \nonumber \\
\tilde{P}_{M-j}(M) &=&
\left(\frac{M^{j-1}}{(2j)!!} + O(M^{j-2})\right)\label{eq:Pol11}
\end{eqnarray}
Here we attach a piece of code in the symbolic computation language Mathematica. 
This code solves the recursion relations (\ref{eq:PolynomialsRecurs}) and verifies the result (\ref{eq:InverseVandermonde}) for the inverse of the Vandermonde matrix.

\begin{verbatim}
M =.;MMAX = 20;MyP[0, M_] = 1;

id = OpenWrite["Polynomials.dat"];

SetOptions[id, FormatType -> TeXForm];

Do[

    MyP[j, M_] = Factor[-Sum[(Sum[ MyP[j - p, l] (-l)^{}p, {p, 1, j}]), {l, 1, M}]];

    Write[id, "P(M-", j, ",M) = ", MyP[j, l] /. {l :> M}, FormatType -> TeXForm];

    Print["P(M-", j, ",M) = ", Simplify[(MyP[j, l] /. {l :> M})/Product[(M - qq), {qq, -1, j - 1}]]];
    , {j, 1, MMAX}];

Close[id];M =.;

MyPP[j_, M_] := MyP[j, l] /. {l :> M};

MyInverse[j_, n_, M_] := (-1)^{}(j + n)/((n - 1)!(M - n)!)Sum[ MyPP[j - p, M](-n)^{}p, {p, 0, j}];

MyVandemonde[M_] := Table[ j^{}(n - 1), {j, 1, M}, {n, 1, M}];

M = MMAX;MatrixForm[MyVandemonde[M]]

MM = Inverse[MyVandemonde[M]];

MatrixForm[MM]

Do[

    Print["Checking M-", j, "th row: ", Table[MyInverse[j, n, M], {n, 1, M}] - MM[[M - j]]];

    , {j, 1, MMAX - 1}];

\end{verbatim}

\appsection{}\label{AppD} 

In this appendix, we present a starting point for determining both necessary and sufficient conditions for the divergent behaviour discussed in the paper.

We solve the system of equations (\ref{eq:SysEq}) by Gaussian elimination.
We eliminate the $(M-p)$\textsuperscript{th} row, for $p=1,\dots,M-s-1$, apply the identity
$x^n-y^n=(x-y)\sum\limits_{t=1}^n x^{n-t} y^t$, divide the equation by $(n-M+p-1)$ and easily arrive at the following result:
\begin{equation}
\left(
\begin{array}{rrrrrrr}
(-1)^M, & \dots,  & \dots,  & \dots, & \dots,  &\sum\limits_{p=j  }^M (-1)^p P_p S_{M-1}^M(p-j ),  & \dots \\
0,      & (-1)^M, & \dots,  & \dots, & \dots,  &\sum\limits_{p=j+1}^M (-1)^p P_p S_{M-2}^M(p-j-1), & \dots \\
0,      & 0,      & (-1)^M, & \dots, & \dots,  &\sum\limits_{p=j+2}^M (-1)^p P_p S_{M-3}^M(p-j-2), & \dots \\
\vdots \\
0,      & 0,      & \dots,  & 0,     & (-1)^M, &\sum\limits_{p=j+M-s-1}^M (-1)^p P_p S_{s}^M(p-j-M+s+1), & \dots 
\end{array}
\right)
\label{eq:GaussElim}
\end{equation}
The rows correspond to $n=M,\dots,s+1$ (from top to bottom) and the columns correspond to $j=M, \dots, 1$ (from left to right). In the bottom row there are $(M-s-1)$ zeros in columns $M, M-1, \dots, s+2$.
For brevity we have dropped the argument $M$ in the $P_p(M)$ polynomials, i.e.\ we have $P_p := P_p(M)$.
The coefficients $S_{s}^M(x)$ read:
\begin{eqnarray}
S_{s}^M(x) &=& \sum\limits_{\sum\limits_{q=1}^{M-s} \theta_q = x} \prod\limits_{q=s+1}^M q^{\theta_{q-s}}
=
\sum\limits_{j=s+1}^{M} j S^j_s(x-1)
\label{eq:CoeffsRecurRels}\\
&=&\frac{(-1)^{M-s-1}}{ (M-s-1)!} \sum\limits_{q=0}^{M-s-1} (-1)^q C^{M-s-1}_q (s+q+1)^{x+M-s-1}
\\
&=& \left(\prod\limits_{q=0}^{x-1} (M-s+q)\right) \left(\sum\limits_{q=0}^x M^{x-q} Q_q^x(s)\right)
\label{eq:CoeffsRecurRels2}\\
&=&
\left\{
\begin{array}{cc}
1   & \mbox{if $s=M$}   \\
M^x & \mbox{if $s=M-1$} \\
M^{x+1} - (M-1)^{x+1} & \mbox{if $s=M-2$} \\
\frac{1}{2} \left(M^{x+2} - 2(M-1)^{x+2} + (M-2)^{x+2}\right) & \mbox{if $s=M-3$} \\
\frac{1}{6} \left(M^{x+3} - 3(M-1)^{x+3} + 3 (M-2)^{x+3} - (M-3)^{x+3}\right) & \mbox{if $s=M-4$} 
\end{array}
\right.
\end{eqnarray}
The quantities $Q_q^x(s)$ are polynomials of order $x$; they satisfy the following recursion relations:
\begin{equation}
Q^{x}_q(s) = \sum\limits_{0\le q_2\le q_1\le q} C^{x-q_1}_{q-q_1} (-s)^{q-q_1} 
A^{x-q_2}_{x-q_1}(s) Q^{x-1}_{q_2}(s)
=
\left\{
\begin{array}{rrrr}
1 \\
\frac{1}{2} & \frac{1+s}{2} &\\
\frac{1}{8} & \frac{6s+7}{24} & \frac{3 s^2 + 5 s+2}{24} \\
\frac{1}{48}& \frac{3s+4}{48} & \frac{3s^2+6s+3}{48} & \frac{s^3+2s^2+s}{48}
\end{array}
\right.
\label{eq:RecurRelsCoeffs}
\end{equation}
for $x=1,2,3$ (row-wise) and $q=0,\dots,x$ (column-wise) subject to $Q^0_0(s) = 1$.
Here the coefficients $A^{x}_{q}(s)$ are defined in (\ref{eq:PolynomCoeffs}).
The result (\ref{eq:RecurRelsCoeffs}) follows from inserting (\ref{eq:CoeffsRecurRels2}) into the second equality on the right hand side in (\ref{eq:CoeffsRecurRels}) and performing the sum over $j$ using (\ref{eq:BinomIdentity2}).

From (\ref{eq:GaussElim}) we see that the solution to the system of equations (\ref{eq:SysEq})
depends on $s$ parameters. We choose the first $s$ reduced probabilities 
$\left(\frac{{\mathbb P}_j - {\mathbb P}_0}{j}\right)_{j=1}^s$ as the parameters.
In each equation in (\ref{eq:SysEq}) we move the last $s$ terms onto the right-hand side and we 
obtain a system of $M-s$ equations with an upper triangular quadratic matrix. We eliminate the consecutive variables
and obtain the following solution:
\begin{eqnarray}
\left(\frac{{\mathbb P}_j - {\mathbb P}_0}{j} \right)_{j=s+1}^M
&=&
\sum\limits_{q=1}^s
\left(
\sum\limits_{\tilde{q}=0}^{j-s-1} (-1)^{\tilde{q}+1}
\sum\limits_{j-1 \ge j_0 > \dots > j_{\tilde{q}-1} \ge s+1}
\frac{\prod\limits_{l=-1}^{\tilde{q}-1}(j_l,j_{l+1})}{\prod\limits_{l=-1}^{\tilde{q}-1} (j_l,j_l)} 
\right)_{j=s+1}^M \nonumber \\
&& \cdot \left(\frac{{\mathbb P}_q - {\mathbb P}_0}{q}\right) \label{eq:ReducedProbabilitiesFinal}\\
&=& 
\sum\limits_{q=1}^s
\vec{{\mathfrak V}}_q
\cdot
\left(\frac{{\mathbb P}_q - {\mathbb P}_0}{q}\right) \label{eq:ReducedProbabilitiesFinal1}
\end{eqnarray}
where 
\begin{equation}
(i,j) := \sum\limits_{p=M+j-i}^M (-1)^p P_p S^M_{i-1}(p-j+i-M)
\label{eq:MatrElems}
\end{equation}
for $j\le i$ and $i=M,\dots,s+1$ and subject to $j_{-1} = j$ and $j_{\tilde{q}} = q$
The step (\ref{eq:ReducedProbabilitiesFinal})$\to$(\ref{eq:ReducedProbabilitiesFinal1}) follows from the fact that
$(i,i) = (-1)^M$.
Here we denoted:
\begin{eqnarray}
\vec{{\mathfrak V}}_q &:=& 
\left(
\sum\limits_{\tilde{q}=0}^{j-s-1} (-1)^{(\tilde{q}+1)(M-1)}
\sum\limits_{j-1 \ge j_0 > \dots > j_{\tilde{q}-1} \ge s+1}
\prod\limits_{l=-1}^{\tilde{q}-1}(j_l,j_{l+1})
\right)_{j=s+1}^M
\label{eq:Mode}
\end{eqnarray}
The sum in (\ref{eq:Mode}) contains $2^{j-q}$ terms.

Now we compute the matrix elements $(i,j)$ in (\ref{eq:MatrElems}) for big values of  $M$. 
We have:
\begin{eqnarray}
(i,j) &=& \sum\limits_{p=M+j-i}^M (-1)^p 
\frac{(M+1)!}{p!}(\frac{M^{M-p-1}}{2^{M-p} (M-p)!} + O(M^{M-p-2}))
\times \nonumber \\
&&\frac{(p-j)!}{(M-i)!}
\left(\sum\limits_{q=0}^{p-j+i-M} M^{p-j+i-M-q} (\frac{C^{p-j+i-M}_q}{(2(p-j+i-M))!!} (i-1)^q + O(i^{q-1}))\right)
\label{eq:LargeMLimitModes}\\
&=& \sum\limits_{p=M+j-i}^M (-1)^p 
\frac{(M+1)!}{p!} \frac{M^{-1-j+i}}{2^{M-p} (M-p)!}  
\left( \frac{(p-j)!}{(M-i)!}\right) \times \nonumber \\
&& \left(\frac{1}{(2(p-j+i-M))!!}\right)
(1 + \frac{i-1}{M})^{p-j+i-M} 
\label{eq:LargeMLimitModes1}\\
&=& 
(M+1) M^{-1-j+i}
\sum\limits_{p=M+j-i}^M (-1)^p
\left( \frac{M!}{p!(M-p)!} \right)
\times \nonumber \\ 
&& \left( \frac{(p-j)!}{(M-i)!(p-j+i-M)!}\right)
\frac{(1 + \frac{i-1}{M})^{p-j+i-M} }{2^{-j+i}}
\label{eq:LargeMLimitModes2}\\
&=&
\frac{M^{i-j}}{2^{-j+i}}
(-1)^{M-i+j}
\sum\limits_{p=0}^{i-j} (-1)^p C^M_{i-j-p} C^{M-i+p}_{M-i} 
(1 + \frac{i-1}{M})^{p}
\label{eq:LargeMLimitModes3}\\
&=&
(\frac{M}{2})^{i-j}
(-1)^{M-i+j} \times
\left\{
\begin{array}{ll}
C^{i-1}_{j-1} & \mbox{if $\lim_{M\rightarrow \infty} \frac{i}{M} = 0$} \\
C^{M}_{i-j} 
F_{2,1}\left[
{\tiny\begin{array}{ll} -i+j & 1-i+M \\ 1-i+j+M \end{array}}
;(1+\gamma)\right] & \mbox{if $\lim_{M\rightarrow \infty} \frac{i}{M} = \gamma > 0$} 
\end{array}
\right.
\label{eq:LargeMLimitModes4}
\end{eqnarray}
In (\ref{eq:LargeMLimitModes}) we used (\ref{eq:CoeffsRecurRels2}) along with (\ref{eq:RecurRelsCoeffs})
to transform the quantities $S^M_{i-1}(p-j+i-M)$
and (\ref{eq:PolynomialsRecurs}) along with (\ref{eq:Pol11}) to transform the quantities $P_p(M)$.
In (\ref{eq:LargeMLimitModes1}) 
we neglected the lower order powers of $M$, we simplified the expression, and we performed the sum over $q=0,\dots,p-j+i-M$ using the binomial expansion formula.
In (\ref{eq:LargeMLimitModes2}) 
we simplified the expression.
In (\ref{eq:LargeMLimitModes3}) we substituted for $p-M-j+i$
and simplified the expression.
Finally in (\ref{eq:LargeMLimitModes3}) we performed the sum over $p$ in the case of big values of $i$
and we expressed that sum through the hypergeometric function otherwise.
Since the latter sum is a hypergeometric sum and since necessary and sufficient conditions for such sums to be expressed in closed form are known (see \cite{Zeilberger} and references therein) it would be interesting to check if a closed form expression can be also found in the latter case.

Now we compute the quantities in (\ref{eq:Mode}) in the large-$M$ limit.
We insert the top result in (\ref{eq:LargeMLimitModes4}) into (\ref{eq:Mode})
and we do the sums over the ordered sequences by using identity (\ref{eq:SumSimplexV}) for $p_1 = s-q$ and $p_2=0$.
The final result reads $\vec{{\mathfrak V}}_q = \left( {\mathfrak V}_{q,j} \right)_{j=s+1}^M$ where:
\begin{equation}
{\mathfrak V}_{q,j} =
(-\frac{M}{2})^{j-q}
\left( 
C^{j-1}_{q-1}  \sum\limits_{l=1}^{j-s} C^{j-s+1}_{l+1} l^{j-q}
+
1_{s>q}
C^{j-1}_{s-1} C^s_{q-1}
\mathop{\sum\limits_{l_1=1,\dots,s-q+1}}_{l_2=1,\dots,j-s+2}
\beta^{j-s,s-q}_{l_1,l_2}
l_1^{s-q} l_2^{j-s} (-1)^{l_1+l_2}
\right)
\label{eq:ModesProbs}
\end{equation}
where
\begin{equation}
\beta^{j-s,s-q}_{l_1,l_2}
=
\sum\limits_{0 < q_1 < j-s- 1} C^{q_1+1}_{l_2} C^{(j-s-q_1)\wedge (s-q+1)}_{l_1+1}
\end{equation}


\end{document}